\documentclass[12pt,a4paper]{article}
\usepackage{amsfonts,amsthm,amsmath,amssymb,amscd}
\usepackage{fancybox}
\usepackage[latin1]{inputenc}
\usepackage[usenames]{color}
\usepackage{graphicx}
\usepackage{enumerate}
\usepackage{float}
\usepackage{fullpage}
\usepackage{epsfig,amssymb} 
\usepackage{subcaption}

\newtheorem{theorem}{Theorem}[section]
\newtheorem{corollary}[theorem]{Corollary}
\newtheorem{lemma}[theorem]{Lemma}
\newtheorem{proposition}[theorem]{Proposition}

\newtheorem{rem}[theorem]{Remark}

\begin{document}
	
\title{\textbf{Closed $p$-Elastic Curves in Spheres of $\mathbb{L}^3$}}

\author{A. P\'ampano, M. Samarakkody and H. Tran}
\date{\today}

\maketitle

\begin{abstract}
For every $p\in\mathbb{R}$, we study $p$-elastic curves in the hyperbolic plane $\mathbb{H}^2$ and in the de Sitter $2$-space $\mathbb{H}_1^2$. We analyze the existence of closed $p$-elastic curves with nonconstant curvature showing that in the hyperbolic plane $\mathbb{H}^2$ these curves exist provided that $p>1$, while in the de Sitter $2$-space $\mathbb{H}_1^2$ the restriction $p<0$ must be satisfied.\\

\end{abstract}

\section{Introduction}

Variational problems for curves whose energy densities depend on the curvature of the curve are ubiquitous in Differential Geometry, Geometric Analysis and Mathematical Physics. Additionally, these type of problems have applications, not only in Mathematics, but also in other fields ranging from Biology and Physics to Architecture and Art.

Their study originated from the pioneering works of the Bernoulli family and L. Euler about the theory of elasticity. Their analysis was essential in the early developments of the Calculus of Variations. In particular, what nowadays are called \emph{$p$-elastic curves} first appeared in a letter from D. Bernoulli to L. Euler (\cite{T}) in 1738. In this correspondence, D. Bernoulli proposed to investigate extrema of the functionals (in modern notation)
$$\mathbf{\Theta}_p(\gamma):=\int_\gamma \kappa^p\,ds\,,$$
where $\kappa$ denotes the curvature of the curve $\gamma$ and $s$ is its arc length parameter. 

Despite their ancient origin, $p$-elastic curves are still a very vital field of research. Even though some particular cases have been well known for a long time, studying them in general gives a rich family of curves which has yet to be fully understood. For instance, the case $p=2$ is the classical bending energy. Its minimizers serve as a simple model to describe the shape of a thin elastic rod with circular cross section and uniform density, naturally straight and prismatic when unstressed and which is being held bent by external forces and moments acting at its ends alone. The planar configurations were classified by L. Euler in 1744 (\cite{E}), although some special cases were already known to J. Bernoulli around 1694 (\cite{B}). More recently, the bending energy has been extended to different ambient spaces (see, for instance, \cite{BG,J,LS}). For a description of some other particular cases and their applications, we refer the reader to the Introduction of \cite{GPT} and the references therein (see also the Introduction of the more recent paper \cite{P0}).

The existence of smooth\footnote{Throughout this paper we will consider smooth curves, although $\mathcal{C}^4$ regularity is enough to deduce our results.} closed $p$-elastic curves with nonconstant curvature in Riemannian and Lorentzian $2$-space forms is an interesting problem which is only partially solved (\cite{ABG,ABG0,GPT,LS0,LS,LP,MOP,MP0,MP}). Except for the trivial variational problem corresponding with $p=1$, $p$-elastic curves in the Euclidean plane $\mathbb{R}^2$ are generically not closed (\cite{LP}). In the Lorentz-Minkowski plane $\mathbb{L}^2$ closed curves must have a change of causal character (\cite{BY}). In the round $2$-sphere $\mathbb{S}^2$ there exist closed $p$-elastic curves with nonconstant curvature if and only if $p=2$ or $p\in(0,1)$ (\cite{AGP,GPT,LS,MOP}). In the latter case, for every pair of relatively prime natural numbers $(n,m)$ satisfying $m<2n<\sqrt{2}\,m$, there exists a closed $p$-elastic curve with nonconstant curvature. The natural numbers $n$ and $m$ represent geometric invariants of the curve: the number $n$ is the winding number, while $m$ describes the symmetry of the curve with respect to the action of a discrete group of rotations. Furthermore, there is a strong numerical evidence that a one-to-one correspondence exists between pairs of relatively prime natural numbers $(n,m)$ satisfying $m<2n<\sqrt{2}\,m$ and closed $p$-elastic curves (for $p\in(0,1)$ fixed) with nonconstant curvature in $\mathbb{S}^2$. This assertion has been analytically shown for the case $p=1/2$ employing well known formulas for the derivatives of elliptic integrals (\cite{AGP}).

Unfortunately, less is known in the remaining $2$-space forms, namely, in the hyperbolic plane $\mathbb{H}^2$ and in the de Sitter $2$-space $\mathbb{H}_1^2$. For instance, in the hyperbolic plane $\mathbb{H}^2$ and for $p>1$ integer, the existence of closed $p$-elastic curves has already been shown (\cite{ABG}). However, only for the cases $p=2$ and $p=3$ there are known formal results regarding the restrictions for their existence (\cite{ABG0,LS0,LS}). The reason why only those two cases have been formally analyzed is that its study can be related to elliptic integrals and functions, while in the general case the theory for the involved functions is not completely developed. Interestingly, in these cases the condition $m<2n<\sqrt{2}\,m$ is again essential. Indeed, closed $p$-elastic curves ($p=2$ and $p=3$) with nonconstant curvature in $\mathbb{H}^2$ are in one-to-one correspondence with the pairs of relatively prime natural numbers $(n,m)$ satisfying that condition (\cite{ABG0,LS0,LS}).

Motivated by the above mentioned results, this work intends to investigate whether or not closed $p$-elastic curves with nonconstant curvature in $\mathbb{H}_\epsilon^2$, $\epsilon=0,1$, exist for a given real number $p\in\mathbb{R}$. In this sense, this paper may be seen as a natural continuation of the works \cite{ABG,GPT}. However, it is not merely an extension of the results and techniques of these works. On one hand, despite obvious formal analogies, the hyperbolic and the pseudo-hyperbolic cases present substantial differences with respect to the spherical case due to the different topologies of their respective isometry groups. In short, the momentum vector field in the hyperbolic and pseudo-hyperbolic cases belongs to the Lorentz-Minkowski space $\mathbb{L}^3$ in contrast to the spherical case, for which the momentum is a vector field in the Euclidean space $\mathbb{R}^3$. As a consequence, the different possibilities for its causal character will play an essential role. Equivalently, in terms of Killing vector fields along curves, both in $\mathbb{H}^2$ and $\mathbb{H}_1^2$ we may encounter rotational, hyperbolic and parabolic types, while in $\mathbb{S}^2$ only the rotational type exists. Additionally, considering real values of $p\in\mathbb{R}$ and not only natural values introduces technical challenges in the analysis of the associated variational problems. For instance, generally the involved integrals and functions will no longer be elliptic. 

Circumventing these difficulties, in this paper, we show that closed $p$-elastic curves with nonconstant curvature in the hyperbolic plane $\mathbb{H}^2$ exist if and only if $p>1$ is any real number. Moreover, we prove that for every pair of relatively prime natural numbers $(n,m)$ satisfying $m<2n<\sqrt{2}\,m$ there exists a closed $p$-elastic curve ($p>1$) with nonconstant curvature. On the other hand, we show that if closed $p$-elastic curves in the de Sitter $2$-space $\mathbb{H}_1^2$ exist, then $p<0$ must hold. We also study the asymptotic behavior of the closure condition. We compute analytically one of the limits obtaining the same result as in the hyperbolic case and complete the argument for the other limit numerically, presenting evidence that for every pair of relatively prime natural numbers $(n,m)$ satisfying $m<2n<\sqrt{2}\,m$ there exists a closed $p$-elastic curve ($p<0$) with nonconstant curvature. In addition, our numerical experiments strongly support the \emph{ansatz} that, for both ambient spaces, there is a one-to-one correspondence between closed $p$-elastic curves (for suitable $p$ fixed) and the pairs $(n,m)$ satisfying $m<2n<\sqrt{2}\,m$.

The findings of the present paper, in combination with the results in the spherical case obtained in \cite{GPT} (with partial contributions of \cite{AGP,MOP}), clarify the problem regarding the existence of closed $p$-elastic curves with nonconstant curvature in Riemannian and Lorentzian $2$-space forms. Indeed, with the remarkable exception of the classical bending energy $\mathbf{\Theta}_2$ acting on curves immersed in $\mathbb{S}^2$, closed $p$-elastic curves with nonconstant curvature are classified as follows: let $(n,m)$ be a pair of relatively prime natural numbers satisfying $m<2n<\sqrt{2}\,m$, then for every $p\in\mathbb{R}\setminus\{0,1\}$ there exists a (space-like) convex closed $p$-elastic curve with nonconstant curvature $\gamma_{n,m}$. Moreover, the curve $\gamma_{n,m}$ is immersed in:
\begin{enumerate}[(i)]
	\item The hyperbolic plane $\mathbb{H}^2$ if $p>1$,
	\item The round $2$-sphere $\mathbb{S}^2$ if $p\in(0,1)$, or
	\item The de Sitter $2$-space $\mathbb{H}_1^2$ if $p<0$.
\end{enumerate}

For the sake of completeness we mention here that the only closed $p$-elastic curve for $p=0$ has constant curvature (in fact, it is a geodesic) and it may be immersed in either $\mathbb{S}^2$ or $\mathbb{H}_1^2$ as the ``equator'' (i.e., in the quadric model for $\mathbb{S}^2$ or $\mathbb{H}_1^2$, it is the parallel located at height $z=0$, up to isometries). The case $p=1$ gives rise to a trivial variational problem and there are no critical curves in $\mathbb{S}^2$, $\mathbb{H}^2$, or $\mathbb{H}_1^2$.

The above classification of convex closed $p$-elastic curves with nonconstant curvature highlights an interesting phenomenon. For fixed pair of relatively prime natural numbers $(n,m)$ satisfying $m<2n<\sqrt{2}\,m$, one can consider the evolution of the corresponding convex closed $p$-elastic curve $\gamma_{n,m}$ as $p$ varies in the real line. Understanding $\mathbb{H}^2$, $\mathbb{S}^2$, and $\mathbb{H}_1^2$ as quadrics in $\mathbb{R}^3$, as $p$ decreases from $+\infty$ to $1$, the curve $\gamma_{n,m}$ goes shrinking while it comes down from the top part of the hyperbolic plane $\mathbb{H}^2$, and as $p\to 1^+$ it converges to the pole $(0,0,1)\in\mathbb{R}^3$. Then, as $p$ decreases from $1$ to $0$, $\gamma_{n,m}$ (which is now immersed in $\mathbb{S}^2$) goes expanding from the north pole $(0,0,1)\in\mathbb{R}^3$ when $p\to 1^-$ to the equator $\mathbb{S}^2\cap\{z=0\}$ as $p\to 0^+$. This equator is the only closed $p$-elastic curve for $p=0$ and it belongs to both $\mathbb{S}^2$ and $\mathbb{H}_1^2$. Finally, when $p<0$ decreases, $\gamma_{n,m}$ expands from the equator of $\mathbb{H}^2_1$ (which coincides with that of $\mathbb{S}^2$) as $p\to 0^-$ while keeps coming down to the bottom part of the de Sitter $2$-space $\mathbb{H}_1^2$. For an illustration of this evolution see Figure \ref{Evolution3} below.

\begin{figure}[h!]
	\centering
	\includegraphics[height=7.5cm]{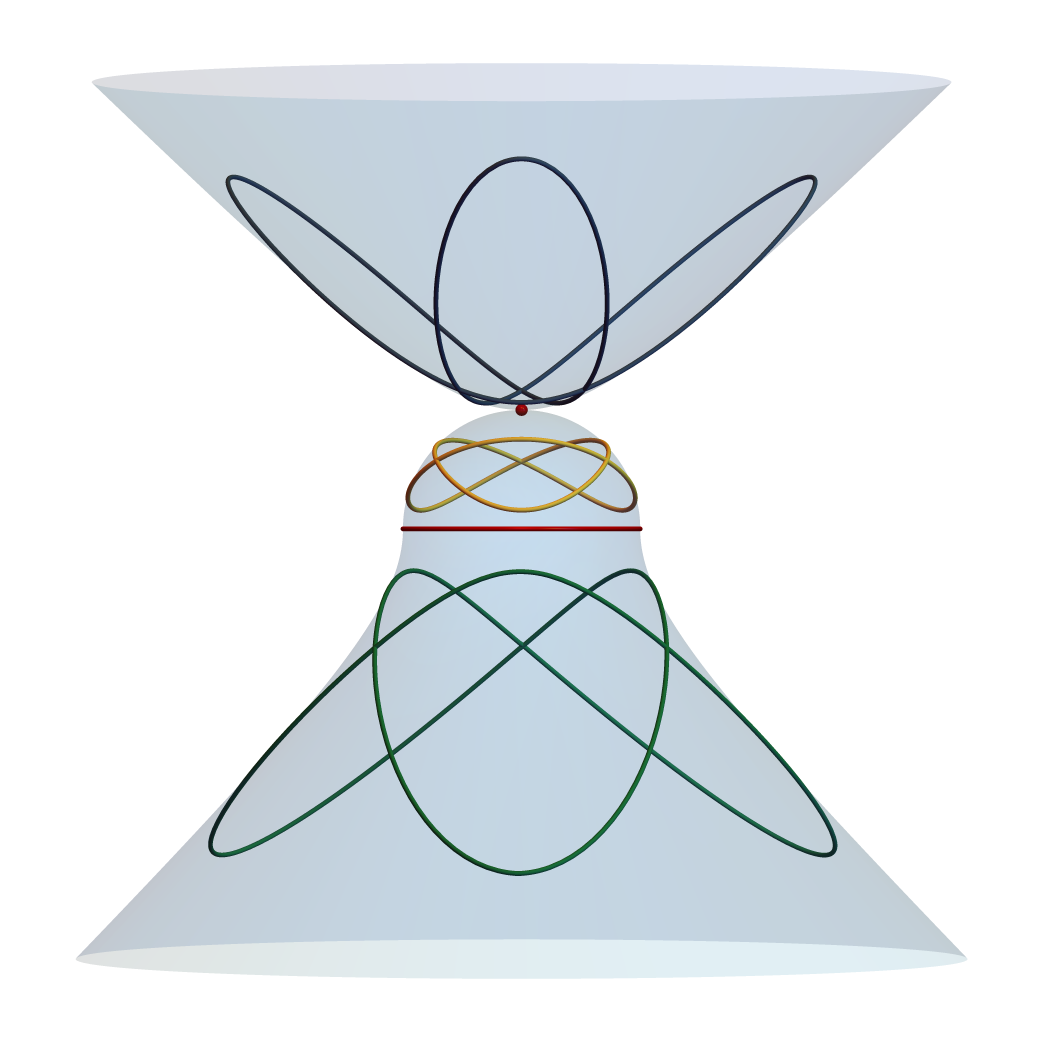}
	\caption{Evolution of closed $p$-elastic curves of type $\gamma_{2,3}$. These curves are represented in the quadric models for $\mathbb{H}^2$, $\mathbb{S}^2$, and $\mathbb{H}_1^2$. In black, the $p$-elastic curve $\gamma_{2,3}$ for $p=2$ immersed in $\mathbb{H}^2$; in yellow, the $p$-elastic curve $\gamma_{2,3}$ for $p=0.2$ immersed in $\mathbb{S}^2$; and, in green, the $p$-elastic curve $\gamma_{2,3}$ for $p=-1$ immersed in $\mathbb{H}_1^2$. The red point is the pole $(0,0,1)\in\mathbb{R}^3$ and the red circle is the equator (i.e., the circle of radius one contained in the plane $\{z=0\}$ and centered at the origin).}
	\label{Evolution3}
\end{figure}

\section{Preliminaries}

Let $(x,y,z)$ be the standard coordinates of $\mathbb{R}^3$. The space $\mathbb{R}^3$ endowed with the Euclidean metric $\overline{g}=dx^2+dy^2+dz^2$ becomes the Euclidean space $\mathbb{E}^3$. The set of points equidistant to the origin of $\mathbb{R}^3$ with the induced metric are the round spheres 
$$\mathbb{S}^2(\rho)=\{(x,y,z)\in\mathbb{R}^3\,\lvert\,\overline{g}\left((x,y,z),(x,y,z)\right)=x^2+y^2+z^2=1/\rho\}\,,$$ 
of constant sectional curvature $\rho>0$.

The \emph{Lorentz-Minkowski $3$-space} $\mathbb{L}^3$ is $\mathbb{R}^3$ endowed with the canonical Lorentzian (i.e., index one) metric
$$g\equiv \langle\cdot,\cdot\rangle=dx^2+dy^2-dz^2\,.$$
In $\mathbb{L}^3$, there are two essentially different immersed non-degenerate surfaces which play the role of the round spheres $\mathbb{S}^2(\rho)$, in the sense that they are the set of points equidistant to the origin with respect to the Lorentzian metric $g$. The induced metric to these surfaces may be Riemannian or Lorentzian. Hence, depending on the causal character of these surfaces we have the hyperbolic planes $\mathbb{H}^2(\rho)$ and the de Sitter 2-spaces $\mathbb{H}_1^2(\rho)$.

The \emph{hyperbolic plane} $\mathbb{H}^2(\rho)$, $\rho<0$, is the space-like surface of $\mathbb{L}^3$ given by the top part of the hyperboloid of two sheets
$$\mathbb{H}^2(\rho)=\{(x,y,z)\in\mathbb{R}^3\,\lvert\,\langle(x,y,z),(x,y,z)\rangle=x^2+y^2-z^2=1/\rho\,,\,z>0\}\,,$$
endowed with the Riemannian metric induced from $\mathbb{L}^3$. The surface $\mathbb{H}^2(\rho)$ has constant sectional curvature $\rho<0$. Without loss of generality, from now on we will assume $\rho=-1$ and use the notation $\mathbb{H}_0^2\equiv\mathbb{H}^2(-1)$. For visualization purposes, we identify $\mathbb{H}_0^2$ with the unit disk $\mathbb{D}\subset\mathbb{R}^2$ centered at the origin and endowed with the Poincar\'e metric by means of the isometry
$$(x,y,z)\in\mathbb{H}_0^2\subset\mathbb{R}^3\longmapsto \frac{1}{1+z}(x,y)\in\mathbb{D}\subset\mathbb{R}^2\,.$$
This identification gives rise to the Poincar\'e disk model for $\mathbb{H}_0^2$.

The \emph{de Sitter $2$-space} $\mathbb{H}_1^2(\rho)$, $\rho>0$, is the time-like surface of $\mathbb{L}^3$ given by the hyperboloid of one sheet
$$\mathbb{H}^2_1(\rho)=\{(x,y,z)\in\mathbb{R}^3\,\lvert\,\langle(x,y,z),(x,y,z)\rangle=x^2+y^2-z^2=1/\rho\,\}\,,$$
endowed with the Lorentzian metric induced from $\mathbb{L}^3$. The surface $\mathbb{H}^2_1(\rho)$ has constant sectional curvature $\rho>0$. We may assume from now on that $\rho=1$ and denote $\mathbb{H}^2_1(1)$ simply by $\mathbb{H}^2_1$. For visualization purposes, we identify the bottom half of $\mathbb{H}_1^2$, namely, $\mathbb{H}_1^2\lvert_-=\mathbb{H}_1^2\cap\{z<0\}$, with the once punctured unit disk $\mathring{\mathbb{D}}=\{(x,y)\in\mathbb{R}^2\,\lvert\,0\neq x^2+y^2<1\}$ via the diffeomorphism
$$(x,y,z)\in\mathbb{H}_1^2\lvert_-\subset\mathbb{R}^3\longmapsto \frac{1}{x^2+y^2}(x,y)\in\mathring{\mathbb{D}}\subset\mathbb{R}^2\,.$$
By considering on $\mathring{\mathbb{D}}$ the pull-back metric, above map becomes an isometry and, hence, $\mathring{\mathbb{D}}$ with that metric is a model for the bottom half of $\mathbb{H}^2_1$.

\begin{rem}
	Although with the above identification we are only representing the bottom half of $\mathbb{H}^2_1$, this will be enough for our purposes (c.f., first property of Proposition \ref{properties}).
\end{rem}

Throughout this paper, we will analyze both the hyperbolic plane $\mathbb{H}_0^2$ and the de Sitter $2$-space $\mathbb{H}_1^2$ in a unified manner. Thus, it is convenient to denote by $\mathbb{H}_\epsilon^2$, $\epsilon=0,1$, the hyperbolic plane $\mathbb{H}^2_0$ and the de Sitter $2$-space $\mathbb{H}_1^2$. With this notation, the sectional curvature of $\mathbb{H}_\epsilon^2$ is $\rho=(-1)^{\epsilon+1}$.

Let $\gamma:\widetilde{I}\subseteq\mathbb{R}\longrightarrow\mathbb{H}_\epsilon^2$ be a smooth immersed curve. If $\epsilon=0$, we say that $\gamma$ is an \emph{hyperbolic curve} and, when $\epsilon=1$, the curve $\gamma$ is said to be a \emph{pseudo-hyperbolic curve}. The velocity vector of $\gamma$ is represented by $\dot{\gamma}(t)=d\gamma/dt(t)$, $t\in\widetilde{I}$. An immersed curve in $\mathbb{H}_\epsilon^2$ is \emph{space-like} (respectively, \emph{time-like} or \emph{null}) if $\langle\dot{\gamma}(t),\dot{\gamma}(t)\rangle>0$ for all $t\in\widetilde{I}$ (respectively, $\langle\dot{\gamma}(t),\dot{\gamma}(t)\rangle<0$ or $\langle\dot{\gamma}(t),\dot{\gamma}(t)\rangle=0$ for all $t\in\widetilde{I}$). 

\begin{rem}\label{curves}
	Throughout this paper we will only consider curves whose causal character does not change as the parameter $t$ moves on the interval $\widetilde{I}\subseteq\mathbb{R}$. Moreover, due to the low dimension of $\mathbb{H}^2_\epsilon$, we will discard null curves since its analysis is trivial.
\end{rem} 

A non-null curve can be parameterized by the natural arc length parameter, which we denote by $s\in I$. Let $\gamma:I\subseteq\mathbb{R}\longrightarrow\mathbb{H}_\epsilon^2$ be a smooth non-null immersed curve parameterized by the arc length parameter $s\in I$. Denote by $T(s):=\gamma'(s)$ the unit tangent vector field along $\gamma(s)$, where $\left(\,\right)'$ denotes the derivative with respect to the arc length parameter $s$, and define the unit normal vector field $N(s)$ along $\gamma(s)$ so that $\{T(s),N(s)\}$ is a positively oriented frame along $\gamma(s)$ of the tangent bundle to $\mathbb{H}_\epsilon^2$, the \emph{Frenet frame}. Since $\gamma$ is non-null, we have that $\langle T(s),T(s)\rangle=\epsilon_1=\pm1$ depending on whether $\gamma$ is space-like ($\epsilon_1=1$) or time-like ($\epsilon_1=-1$). Similarly, $\langle N(s),N(s)\rangle=\epsilon_2=\pm1$.

\begin{rem}
	Since the hyperbolic plane $\mathbb{H}^2_0$ is a space-like surface of $\mathbb{L}^3$ (in other words, the induced metric is Riemannian), hyperbolic curves are always space-like $(\epsilon_1=1)$. Moreover, $\epsilon_2=1$ also holds. On the other hand, the de Sitter $2$-space $\mathbb{H}_1^2$ is a time-like surface and, hence, for pseudo-hyperbolic curves we have $\epsilon_1\epsilon_2=-1$.
\end{rem}

In this setting, the (signed) geodesic \emph{curvature} $\kappa(s)$ of $\gamma(s)$ is defined by the Frenet-Serret equations
\begin{eqnarray}
	\nabla_T T(s)&=&\epsilon_2\kappa(s)N(s)\,,\label{FSE1}\\
	\nabla_T N(s)&=&-\epsilon_1\kappa(s)T(s)\,,\label{FSE2}
\end{eqnarray}
where $\nabla$ denotes the Levi-Civita connection on the tangent bundle to $\mathbb{H}_\epsilon^2$. From the Fundamental Theorem of Curves, the curvature $\kappa(s)$ and the causal characters of the Frenet frame $\{T(s),N(s)\}$ ($\epsilon_1$ and $\epsilon_2$, respectively) completely determine the curve $\gamma(s)$, up to rigid motions of $\mathbb{H}_\epsilon^2$.

If $\kappa(s_o)=0$ for some value of $s_o\in I\subseteq\mathbb{R}$, we say that $\gamma(s_o)$ is an \emph{inflection point}. A curve with, at least, one inflection point will be referred as to \emph{inflectional curve}. Curves with no inflection points have a preferred orientation such that $\kappa>0$ everywhere. In other words, with this choice of orientation, they are \emph{convex curves}. We will always assume that curves with no inflection points are oriented in this preferred way so that they are convex curves.

\section{Critical Curves}

Let $\mathcal{C}^\infty(I,\mathbb{H}_\epsilon^2)$ be the space of smooth non-null immersed curves $\gamma:I\subseteq\mathbb{R}\longrightarrow\mathbb{H}_\epsilon^2$ parameterized by the arc length $s\in I$ and denote by $\mathcal{C}^\infty_*(I,\mathbb{H}_\epsilon^2)$ the subspace of convex curves, i.e., those curves whose curvature is strictly positive $\kappa(s)>0$ for all $s\in I$ (curves with no inflection points suitably oriented). 

For any fixed real number $p\in\mathbb{R}$, the \emph{$p$-elastic functional} is given by
\begin{equation}\label{energy}
	\mathbf{\Theta}_p(\gamma):=\int_\gamma \kappa^p\,ds\,,
\end{equation}
and it acts on $\mathcal{C}^\infty(I,\mathbb{H}_\epsilon^2)$. If $p\in\mathbb{R}\setminus\mathbb{N}$ (we consider here $0\in\mathbb{N}$), $\mathbf{\Theta}_p$ is understood to act on $\mathcal{C}^\infty_*(I,\mathbb{H}_\epsilon^2)$.

Using compactly supported smooth variations and standard arguments involving integration by parts, it can be shown that, regardless of the boundary conditions, critical points for $\mathbf{\Theta}_p$ must satisfy the Euler-Lagrange equation
\begin{equation}\label{EL}
	p\,\frac{d^2}{ds^2}\left(\kappa^{p-1}\right)+\epsilon_1\epsilon_2(p-1)\kappa^{p+1}-\epsilon_2 p \kappa^{p-1}=0\,,
\end{equation}
where $\epsilon_1$ and $\epsilon_2$ are the causal characters of the Frenet frame associated to the critical point. Recall that $\epsilon_1\epsilon_2=(-1)^\epsilon$. Those curves immersed in $\mathbb{H}_\epsilon^2$ whose curvature is a solution of \eqref{EL} are called (free) \emph{$p$-elastic curves}.

\begin{rem}\label{geodesics}
	When $p\in\mathbb{N}$, geodesics, i.e., curves immersed in $\mathbb{H}_\epsilon^2$ satisfying $\kappa=0$ everywhere, are always solutions of \eqref{EL} and, hence, $p$-elastic curves. For later use, we point out here that the only closed geodesic of $\mathbb{H}_\epsilon^2$ arises when $\epsilon=1$, that is, in the de Sitter $2$-space $\mathbb{H}_1^2$ and it is the parallel circle located at height $z=0$. 
\end{rem}

Before proceeding further, we now discuss all the solutions of \eqref{EL} for two special values of $p\in\mathbb{R}$. If $p=0$, the Euler-Lagrange equation \eqref{EL} reduces to $\kappa=0$ and, hence, the only critical curves are geodesics. This is coherent since the functional $\mathbf{\Theta}_0$ is just the length functional. On the other hand, if $p=1$ we deduce from \eqref{EL} that critical curves cannot exist, since $\epsilon_2=\pm1\neq 0$. From now on, we will assume $p\neq 0,1$ and we will analyze nonzero solutions of \eqref{EL}.

We first consider constant solutions of \eqref{EL}. Since we are mainly interested on obtaining closed $p$-elastic curves, we will study the existence of $p$-elastic circles, i.e., circles of $\mathbb{L}^3$ immersed in $\mathbb{H}_\epsilon^2$ whose constant curvature is a solution of \eqref{EL}.

\begin{proposition}\label{constant}
	Let $\gamma$ be a non-geodesic $p$-elastic circle immersed in $\mathbb{H}_\epsilon^2\subset\mathbb{L}^3$. Then, $\gamma$ is space-like and its constant curvature is given by
	\begin{equation}\label{constantcurvature}
		\kappa=\sqrt{\frac{p}{p-1}}\,.
	\end{equation}
	Equivalently, the radius of $\gamma$, viewed as a curve in $\mathbb{L}^3$, is $r=\sqrt{(-1)^\epsilon(p-1)}$. Moreover:
	\begin{enumerate}[(i)]
		\item If $\gamma\subset\mathbb{H}_0^2$ is a hyperbolic curve, then $p>1$ holds.
		\item If $\gamma\subset\mathbb{H}_1^2$ is a pseudo-hyperbolic curve, then $p<0$ holds.
	\end{enumerate}
\end{proposition}
\textit{Proof.} Let $\gamma$ be a circle immersed in $\mathbb{H}_\epsilon^2\subset\mathbb{L}^3$. Notice that circles of $\mathbb{L}^3$ must be contained in space-like planes or, otherwise, their causal character would change as the circle closes (recall that we are only considering curves whose causal character does not change, see Remark \ref{curves}). This shows that $\gamma$ is space-like ($\epsilon_1=1$). In addition, it also implies that the covariant derivative $\overline{\nabla}_TT$ is space-like. Here, we are denoting by $\overline{\nabla}$ the Levi-Civita connection on $\mathbb{L}^3$ (do not confuse with $\nabla$, the Levi-Civita connection on $\mathbb{H}_\epsilon^2$).

Assume that $\gamma$ is a $p$-elastic circle in $\mathbb{H}_\epsilon^2$. Then, its constant curvature $\kappa$ is a solution of \eqref{EL}. Since $\kappa$ is constant, \eqref{EL} reduces to
$$\kappa^{p-1}\left(\left(p-1\right)\kappa^2-p\right)=0\,.$$
Hence,
$$\kappa^2=\frac{p}{p-1}\,,$$
must hold because $\gamma$ is a non-geodesic curve.

To compute the radius of the circle $\gamma$, viewed as a curve in $\mathbb{L}^3$, we use the Gauss Formula to deduce
$$\overline{\kappa}^2=\langle\overline{\nabla}_TT,\overline{\nabla}_TT\rangle=\epsilon_2\kappa^2+\rho=\frac{\epsilon_2p}{p-1}-\epsilon_2=\frac{\epsilon_2}{p-1}\,,$$
where $\overline{\kappa}$ denotes the curvature of $\gamma$ in $\mathbb{L}^3$. Consequently, using the relation $r=1/\overline{\kappa}$ for the radius of curvature of a planar curve, the radii of $p$-elastic circles are given by
$$r=\sqrt{\epsilon_2(p-1)}=\sqrt{(-1)^\epsilon(p-1)}\,,$$
since $\epsilon_2=(-1)^\epsilon$ holds for space-like curves.

Finally, if $\gamma$ is a hyperbolic curve ($\epsilon=0$), then $r=\sqrt{p-1}$ and so $p>1$ must hold. On the other hand, if $\gamma$ is a pseudo-hyperbolic curve ($\epsilon=1$), we have $r=\sqrt{1-p}$ and, hence, $p<1$. In addition, it follows from \eqref{constantcurvature} that $p<0$ must hold. \hfill$\square$

\begin{rem}\label{parallels}
	An alternative characterization of non-geodesic $p$-elastic circles immersed in $\mathbb{H}_\epsilon^2$ is that they are the parallel circles located at heights $z^2=(-1)^\epsilon p$.
\end{rem}

In what follows, we will obtain a conservation law for the Euler-Lagrange equation \eqref{EL}. Since in Proposition \ref{constant} we have described the $p$-elastic circles, we will now focus our attention on nonconstant solutions of \eqref{EL}.

\begin{proposition}\label{firstintegral}
	Let $\gamma(s)$ be an arc length parameterized $p$-elastic curve in $\mathbb{H}_\epsilon^2$ with nonconstant curvature $\kappa(s)$. Then, $\kappa(s)$ is a solution of the first order ordinary differential equation,
	\begin{equation}\label{fi}
		p^2\left(p-1\right)^2\kappa^{2(p-2)}\left(\kappa'\right)^2+\epsilon_1\epsilon_2\left(p-1\right)^2\kappa^{2p}-\epsilon_2p^2\kappa^{2(p-1)}=a\,,
	\end{equation}
	where $a\in\mathbb{R}$ is a constant of integration and $\epsilon_1,\epsilon_2$ are the causal characters of the Frenet frame $\{T(s),N(s)\}$ along $\gamma(s)$.
\end{proposition}
\textit{Proof.} Let $\gamma(s)$ be a $p$-elastic curve with nonconstant curvature $\kappa(s)$. Since the curvature is nonconstant, so is $\kappa^{p-1}(s)$. Hence, its derivative with respect to the arc length parameter $s$ is not identically zero. We can then multiply \eqref{EL} by this derivative, obtaining
$$p\,\frac{d^2}{ds^2}\left(\kappa^{p-1}\right)\frac{d}{ds}\left(\kappa^{p-1}\right)+\epsilon_1\epsilon_2(p-1)\frac{d}{ds}\left(\kappa^{p-1}\right)\kappa^{p+1}-\epsilon_2p\,\frac{d}{ds}\left(\kappa^{p-1}\right)\kappa^{p-1}=0\,.$$
The first and last terms above are exact derivatives, while the second one can be manipulated to get
$$\frac{d}{ds}\left(\kappa^{p-1}\right)=(p-1)\kappa^{p-2}\kappa'\,.$$
Consequently, integrating and multiplying by $2p$ we deduce that
$$p^2\left(\frac{d}{ds}\left(\kappa^{p-1}\right)\right)^2+\epsilon_1\epsilon_2(p-1)^2\kappa^{2p}-\epsilon_2p^2\kappa^{2(p-1)}=a\,,$$
for some real constant of integration $a\in\mathbb{R}$. Using this and expanding the first term in above integral expression we conclude with the result. \hfill$\square$

\begin{rem}\label{momentum}
	A geometric interpretation of the conservation law \eqref{fi} is that the momentum,
	$${\vec \xi}:=p\,\kappa^{p-1}\gamma+p\,\frac{d}{ds}\left(\kappa^{p-1}\right)\gamma'+(1-p)\kappa^p \gamma\times\gamma'\in\mathbb{L}^3\,,$$
	is constant along $p$-elastic curves. More precisely, \eqref{fi} is equivalent to $\langle {\vec \xi},{\vec \xi}\rangle=\epsilon_1a$. In the definition of ${\vec \xi}$, $\times$ denotes the usual vector cross product. An alternative description of the conservation law \eqref{fi} can be given in terms of Killing vector fields along $p$-elastic curves (for more details, see the next section).
\end{rem}

For fixed $p\in\mathbb{R}$, solutions of the conservation law \eqref{fi} depend on two real parameters, namely, the constant of integration $a\in\mathbb{R}$ and another constant of integration that will arise when solving \eqref{fi}, as well as in the causal character of the Frenet frame (which, in particular, encodes the choice of ambient space). By simply shifting the origin of the arc length parameter $s$, the integration constant coming from solving \eqref{fi} can be assumed to be zero, hence, solutions of \eqref{fi} depend, essentially, on just the parameter $a\in\mathbb{R}$ and the causal characters. Consequently, for fixed causal characters, there is a (real) one-parameter family of functions $\kappa_a(s)$ solving \eqref{fi}. As mentioned above, the curvature functions $\kappa_a(s)$ (together with the causal character of the Frenet frame) uniquely determine the $p$-elastic curves of $\mathbb{H}_\epsilon^2$, up to isometries. For the sake of clarity, we will denote by $\gamma_a$ the $p$-elastic curve whose curvature $\kappa_a$ is a solution of \eqref{fi} for the fixed parameter $a\in\mathbb{R}$.

We finish this section by proving that, but for one remarkable exception, all $p$-elastic curves with nonconstant curvature are convex. That is, they have no inflection points.

\begin{proposition}\label{inflectional}
	Let $\gamma_a:I\subseteq\mathbb{R}\longrightarrow\mathbb{H}_\epsilon^2$ be a $p$-elastic curve with nonconstant curvature $\kappa\equiv\kappa_a$ and assume that $\gamma_a$ has an inflection point, i.e., there exists $s_o\in I$ such that $\kappa_a(s_o)=0$. Then, $p=2$ and $a>0$ holds. Moreover, if in addition $\gamma_a$ is complete ($I=\mathbb{R}$) and $\kappa_a$ is periodic, then $\gamma_a$ is a hyperbolic curve (i.e., $\epsilon=0$).
\end{proposition}
\textit{Proof.} Assume that the $p$-elastic curve $\gamma_a$ has, at least, one inflection point. It then follows that $p\in\mathbb{N}$, since for $p\in\mathbb{R}\setminus\mathbb{N}$ all our curves are requested to be convex, i.e., $\kappa>0$ holds. In addition, since $\gamma_a$ has nonconstant curvature, $p\geq 2$ must hold. For the cases $p=0,1$ we have already discussed that there are no $p$-elastic curves with nonconstant curvature.

We next show that if $p>2$ is a natural number, solutions of \eqref{fi} must have $\kappa>0$. By contradiction, suppose that there exists $s_o\in I$, where $I$ is the interval of definition of $\gamma$, such that $\kappa(s_o)=0$. Then, evaluating \eqref{fi} at $s_o$ we conclude that $a=0$ must hold. In this case, we can rewrite \eqref{fi} as
$$\left(\kappa'\right)^2=\frac{\kappa^2}{p^2(p-1)^2}\left(\epsilon_2 p^2-\epsilon_1\epsilon_2(p-1)^2\kappa^2\right).$$
Using $\kappa(s_o)=0$ as the initial condition, we get that the only solution to this initial value problem is $\kappa\equiv 0$, reaching a contradiction. Consequently, $\kappa>0$ must hold.

Let $p=2$. Then, the conservation law \eqref{fi} reduces to
\begin{equation}\label{p2}
	\left(\kappa'\right)^2=\frac{a}{4}+\epsilon_2\kappa^2-\epsilon_1\epsilon_2\frac{1}{4}\kappa^4=\widehat{Q}_a(\kappa)\,.
\end{equation}
Evaluating this at the inflection point $s_o\in I$, we obtain
$$\left(\kappa'(s_o)\right)^2=\frac{a}{4}\,,$$
and, hence, $a\geq 0$ must hold. If $a=0$, an argument as above will prove that $\kappa\equiv 0$, reaching a contradiction. Consequently, we deduce that $a>0$ must hold.

Assume, in addition, that the curvature $\kappa$ is a periodic function. Then, the polynomial $\widehat{Q}_a(\kappa)$ of degree four defined in \eqref{p2} must have, at least, two roots. Clearly, if $\alpha>0$ is one such a root, by symmetry $-\alpha<0$ is another one. Moreover, $\widehat{Q}_a(\kappa)>0$ must hold when $\kappa\in(-\alpha,\alpha)$. A simple analysis of this polynomial then shows that this occurs when $\epsilon_1=1$ (that is the curve is space-like) and either $\epsilon_2=1$ and the curve is hyperbolic or $\epsilon_2=-1$ and $0<a\leq 4$.

We also assume now that the $2$-elastic curve $\gamma_a$ is defined on the entire real line $\mathbb{R}$. Under this extra assumption we will discard the possibility that the curve is pseudo-hyperbolic. Let $\gamma_a$ be a pseudo-hyperbolic space-like curve (i.e., $\epsilon_1=1$ and $\epsilon_2=-1$ hold). From \cite{LS} (see also \cite{NEW}), the vector field (c.f., \eqref{J} below)
$$\mathcal{J}=\kappa^2 T+2\kappa' N\,,$$
along $\gamma_a$ is a Killing vector field along the curve. Therefore, it can be uniquely extended to a classical Killing vector field in $\mathbb{H}_1^2$ (see, for instance, \cite{NEW,LS}). The integral curves of Killing vector fields may be circles, hyperbolas or parabolas and so, in particular, the vector field $\mathcal{J}$ cannot have a change of causal character. We then compute,
$$\langle \mathcal{J},\mathcal{J}\rangle=\kappa^4-4\left(\kappa'\right)^2=4\kappa^2-a\,,$$
where we have used that \eqref{p2} holds. Beginning at the inflection point $\gamma_a(s_o)$, $s_o\in I$, we see that $\mathcal{J}$ is time-like since $\kappa(s_o)=0$ holds and $0<a\leq 4$. The vector field $\mathcal{J}$ becomes null whenever $\kappa=\pm\sqrt{a}/2$ (for simplicity, we will argue for the first value, $s_*\in I$, such that $\kappa(s_*)=\sqrt{a}/2$). It turns out that,
$$\kappa(s_*)=\frac{\sqrt{a}}{2}<\sqrt{2-\sqrt{4-a}}=\alpha\,,$$
that is $\mathcal{J}$ changes causal character before the maximum curvature $\alpha$ is attained and so the values of $s\in I$ at which the maximum curvature is attained cannot be included in the maximal interval of definition of $\gamma_a$. Consequently, $\gamma_a$ cannot be defined on the entire real line $\mathbb{R}$, contradicting our assumption. \hfill$\square$

\begin{rem}
	Proposition \ref{inflectional} shows that inflectional $p$-elastic curves are elastic curves in the hyperbolic plane $\mathbb{H}_0^2$ in the classical sense, that is, $p=2$ and $\epsilon=0$ hold. Furthermore, inflectional elastic curves with nonconstant periodic curvature are (space-like) curves for which $a>0$ holds. This means that their momenta ${\vec \xi}\in\mathbb{L}^3$ are space-like and, hence, these curves are not closed (c.f., \cite{LS}).
\end{rem}

As stated in the previous remark (for a proof see \cite{LS}), inflectional $p$-elastic curves are not closed. Consequently, we will discard them and, from now on, just work with convex $p$-elastic curves.

\section{Critical Curves with Periodic Curvatures}

In order to find closed $p$-elastic curves other than the circles obtained in Proposition \ref{constant}, we first need to check whether or not $p$-elastic curves with nonconstant curvature exist. The periodicity of the curvature is a necessary, but not sufficient (see next section), condition for a curve to be closed.

In this section, we will show that periodic solutions $\kappa_a(s)$ of \eqref{fi} exist provided that suitable restrictions on the values of the constant of integration $a\in\mathbb{R}$ and on the parameter $p\in\mathbb{R}$ are imposed. Then, we will explicitly obtain an arc length parameterization of the associated $p$-elastic curves $\gamma_a(s)$ in terms of just one quadrature and, based on this, describe some of their geometric properties. In particular, we will describe a closure condition.

\subsection{Restrictions and Existence}

We begin describing the necessary restrictions on the parameters $a\in\mathbb{R}$ and $p\in\mathbb{R}$ so that the functions $\kappa_a(s)$ solving \eqref{fi} are periodic. In particular, we have already discussed that for $p=0,1$, there are no nonconstant solutions (periodic or not) to \eqref{EL}.

\begin{theorem}\label{periodic} 
	Let $\gamma_a$ be a convex $p$-elastic curve in $\mathbb{H}_\epsilon^2$ with nonconstant periodic curvature. Then, $\gamma_a$ is a space-like curve ($\epsilon_1=1$), \begin{equation}\label{a*}
		0>a>a_*:=-\left((-1)^\epsilon p\right)^p\left((-1)^\epsilon (p-1)\right)^{1-p}\,,
	\end{equation}
	and:
	\begin{enumerate}[(i)]
		\item If $\gamma_a\subset\mathbb{H}_0^2$ is a hyperbolic curve, then $p>1$ holds.
		\item If $\gamma_a\subset\mathbb{H}_1^2$ is a pseudo-hyperbolic curve, then $p<0$ holds.
	\end{enumerate}
	In both cases, the momentum ${\vec \xi}\in\mathbb{L}^3$ is a time-like vector, i.e., $\langle{\vec \xi},{\vec \xi}\rangle=a<0$.
\end{theorem}
\textit{Proof.} Consider a convex $p$-elastic curve $\gamma_a$ with nonconstant periodic curvature $\kappa_a$. Throughout the proof we will avoid explicitly writing the subindex $a$. Then, $\kappa\equiv\kappa_a>0$ is a periodic solution of the conservation law \eqref{fi}. Due to the different powers involved in the expression of this differential equation, we will distinguish between the cases $p>1$ and $p<1$ (for $p=1$ there are no solutions of \eqref{EL} and so neither of its first integral).

Assume first that $p<1$ holds. In this case, we rewrite \eqref{fi} as
$$\left(\kappa'\right)^2=\frac{\kappa^2}{p^2(p-1)^2}\left(a\kappa^{2(1-p)}-\epsilon_1\epsilon_2(p-1)^2\kappa^2+\epsilon_2p^2\right)=\frac{\kappa^2}{p^2(p-1)^2}\,Q_{p,a}(\kappa)\,.$$
Since $\kappa$ is a periodic solution of the above equation, then $Q_{p,a}(\kappa)=0$ must have two positive solutions and $Q_{p,a}(\kappa)>0$ must hold for all the values of $\kappa$ between those two solutions. Differentiating $Q_{p,a}(\kappa)$ we obtain that $a\neq 0$ or, otherwise, $Q_{p,a}$ would be monotonic for positive values of $\kappa$, which is not compatible with having two positive solutions of $Q_{p,a}(\kappa)=0$. In the case $a\neq 0$, $Q_{p,a}$ may have, at most, one critical point $\kappa_*>0$ which is given by the relation
\begin{equation}\label{k*}
	\kappa_*^{-2p}=\epsilon_1\epsilon_2\frac{1-p}{a}\,.
\end{equation}
Since $Q_{p,a}$ has at most one critical point we need $\lim_{\kappa\to 0}Q_{p,a}(\kappa)=\epsilon_2p^2<0$ and that at $\kappa_*$ a local maximum is attained. From the first condition $\epsilon_2=-1$ and, hence, $\epsilon_1=1$ and we are in the de Sitter $2$-space $\mathbb{H}_1^2$. Combining the second condition with
$$Q''_{p,a}(\kappa_*)=4p(1-p)^2\,,$$
we deduce that $p<0$ must hold in order $Q_{p,a}$ to attain a local maximum at $\kappa_*$. This implies from \eqref{k*} that $a<0$. In addition, $Q_{p,a}$ evaluated at $\kappa_*$ must be positive. Requesting this, we get
$$Q_{p,a}(\kappa_*)=-p\left(\frac{(1-p)^{(p-1)/p}}{(-a)^{-1/p}}+p\right)>0\,,$$
or, since $p<0$, equivalently $a>-(-p)^p(1-p)^{1-p}=a_*$. 

Assume now that $p>1$ holds. For these values of $p\in\mathbb{R}$, we may rewrite \eqref{fi} as
$$\left(\kappa'\right)^2=\frac{1}{p^2(p-1)^2\kappa^{2(p-2)}}\left(a-\epsilon_1\epsilon_2(p-1)^2\kappa^{2p}+\epsilon_2 p^2\kappa^{2(p-1)}\right)=\frac{1}{p^2(p-1)^2\kappa^{2(p-2)}}\,\widetilde{Q}_{p,a}(\kappa)\,.$$
Since $\kappa>0$ is a periodic solution, we have that $\widetilde{Q}_{p,a}(\kappa)=0$ has two positive solutions and for the values of $\kappa$ between them, $\widetilde{Q}_{p,a}(\kappa)>0$ holds. Differentiating $\widetilde{Q}_{p,a}$ we get that it has a unique positive critical point $\kappa_*$, which satisfies
$$\kappa_*^2=\frac{p}{p-1}\,,$$
provided that $\epsilon_1=1$. To the contrary, if $\epsilon_1=-1$, $\widetilde{Q}_{p,a}(\kappa)$ would be monotonic for positive values of $\kappa$ and periodic solutions of \eqref{fi} would not exist. In addition, since for $\epsilon_1=1$ there is just one positive critical point, namely, $\kappa_*$, the function $\widetilde{Q}_{p,a}$ must attain a local maximum there and $\lim_{\kappa\to 0}\widetilde{Q}_{p,a}(\kappa)=a<0$, necessarily. Differentiating $\widetilde{Q}_{p,a}$ twice we conclude that at $\kappa_*$ a local maximum is attained if and only if
$$\widetilde{Q}''_{p,a}(\kappa_*)=-4\epsilon_2\frac{p^p}{(p-1)^{p-3}}<0\,,$$
that is, $\epsilon_2=1$ and we are in the hyperbolic plane $\mathbb{H}_0^2$. Finally, we also have $\widetilde{Q}_{p,a}(\kappa_*)>0$. This gives
$$\widetilde{Q}_{p,a}(\kappa_*)=a+\frac{p^p}{(p-1)^{p-1}}>0\,,$$
which implies $a>-p^p(p-1)^{1-p}=a_*$. This finishes the proof. \hfill$\square$
\\

The converse of Theorem \ref{periodic} also holds. Assuming the corresponding restrictions on the parameters $a,p\in\mathbb{R}$, the curvature $\kappa_a$ of a space-like curve $\gamma_a$ in $\mathbb{H}_\epsilon^2$ solution of \eqref{fi} is a periodic function. This shows the existence of $p$-elastic curves with periodic curvatures. Furthermore, it shows that, under the restrictions of Theorem \ref{periodic}, all solutions of \eqref{fi} are periodic functions, provided that they are defined on their maximal domains.

\begin{proposition}\label{converse} Assume that $p>1$ holds when $\mathbb{H}_\epsilon^2=\mathbb{H}_0^2$ is the hyperbolic plane and that $p<0$ holds when $\mathbb{H}_\epsilon^2=\mathbb{H}_1^2$ is the de Sitter $2$-space. Then, for every fixed real number $a\in\mathbb{R}$ such that $0>a>a_*$ holds, there exists a space-like convex $p$-elastic curve $\gamma_a:I\subseteq\mathbb{R}\longrightarrow\mathbb{H}_\epsilon^2$ with nonconstant curvature $\kappa_a$. If, in addition, $\gamma_a$ is defined on its maximal domain, then it is complete $(I=\mathbb{R})$ and its curvature $\kappa_a$ is a periodic function.
\end{proposition}
\textit{Proof.} Fix the numbers $\epsilon=0,1$, $p\in\mathbb{R}$ and $a\in\mathbb{R}$ such that they satisfy the restrictions imposed in the statement. Fix also $\epsilon_1=1$. For these constant values, it was shown in Theorem \ref{periodic} that there exists two positive values, namely, $\alpha>\beta>0$ such that $Q_{p,a}(\alpha)=Q_{p,a}(\beta)=0$ (respectively, for $\widetilde{Q}_{p,a}(\kappa)$). Pick up one of them to be the initial condition for the first order ordinary differential equation \eqref{fi}, say $\kappa_a(s_o)=\alpha$ with $s_o\in I$. Then, the local existence of solution $\kappa_a(s)$ follows from the standard theory of differential equations. Since the solution satisfies the conservation law \eqref{fi} it also satisfies the Euler-Lagrange equation \eqref{EL}. 

Now, we have the solution $\kappa_a(s)$, the causal characters $\epsilon_1=1$ and $\epsilon_2=(-1)^\epsilon$ of the Frenet frame and the arc length parameter $s\in I$ fixed. The Fundamental Theorem of Curves then guarantees the existence of a unique (up to isometries) curve in $\mathbb{H}_\epsilon^2$. We denote it by $\gamma_a$. The curve $\gamma_a$ is a space-like $p$-elastic curve by construction, since its curvature satisfies the Euler-Lagrange equation.

We next check that the solution $\kappa_a(s)$ is not constant. Assume by contradiction that $\kappa_a(s)\equiv\alpha$ is a constant solution of \eqref{fi}. Then,
$$\epsilon_2\alpha^{2(p-1)}\left((p-1)^2\alpha^2-p^2\right)=a<0\,,$$
holds. However, we have seen in Proposition \ref{constant} that $p$-elastic curves with constant curvature satisfy \eqref{constantcurvature} or they are geodesics. In either case, we reach a contradiction, so $\kappa_a(s)$ is nonconstant. Furthermore, $\kappa_a(s)>0$ for every $s\in I$. Otherwise, according to Proposition \ref{inflectional}, $a>0$ which is not our case. This shows that $\gamma_a(s)$ is convex.

Let $\gamma_a$ be a space-like convex $p$-elastic curve in $\mathbb{H}_\epsilon^2$ with nonconstant curvature $\kappa_a$, where $a\in\mathbb{R}$ is a fixed number such that $0>a>a_*$ holds. The existence of such a curve has just been proved. Moreover, the curvature function $\kappa_a$ must be a solution of the conservation law \eqref{fi}. With the new variables $u=\kappa_a$ and $v=u'=\kappa'_a$, \eqref{fi} becomes
\begin{equation}\label{uv-curve}
	v^2=\frac{1}{p^2(p-1)^2}\left(au^{2(2-p)}-\epsilon_2(p-1)^2u^4-\epsilon_2p^2u^2\right).
\end{equation}
This is a curve in the phase plane representing the orbits of the differential equation \eqref{fi}. Observe that what is inside the parenthesis is, possibly after taking out a common factor, either $Q_{p,a}(u)$ or $\widetilde{Q}_{p,a}(u)$ (c.f., Theorem \ref{periodic}). Indeed, repeating the analysis of these functions as already done in the proof of Theorem \ref{periodic}, we conclude that under the assumptions of the statement, the curve \eqref{uv-curve} is closed.

Furthermore, the curve $C(s)=\left(u(s),v(s)\right)$, $s\in I$, is included in the trace of the closed curve \eqref{uv-curve} and is the bounded integral curve of the smooth vector field
$${\vec X}(u,v)=\left(v,\frac{1}{u}\left((2-p)v^2-\epsilon_2\frac{1}{p}u^4+\epsilon_2\frac{1}{1-p}u^2\right)\right),$$
which is defined on $\{(u,v)\in\mathbb{R}^2\,\lvert\,u>0\}$. Hence, $C(s)$ is smooth and defined on the entire real line $\mathbb{R}$. In particular, so is $u(s)=\kappa_a(s)$ and its associated curve $\gamma_a(s)$. In other words, $I=\mathbb{R}$ is the maximal domain of definition of $\gamma_a$.

In addition, since the vector field ${\vec X}$ has no zeros along $C(s)$, we conclude from the Poincar\'e-Bendixson Theorem that $C(s)$ is itself periodic. This means that $\kappa_a(s)$ is a periodic function. \hfill$\square$
\\

In what follows we will always assume that our curves are defined on their maximal domain. Hence, previous restrictions on the energy parameter $p\in\mathbb{R}$ and on the constant of integration $a\in\mathbb{R}$ would guarantee that $p$-elastic curves have periodic curvatures.

\subsection{Geometric Properties}

We will now obtain a parameterization of convex $p$-elastic curves with nonconstant periodic curvature. Recall that as a consequence of the periodicity $0>a>a_*$ holds (which, for the sake of simplicity, we will omit if unnecessary) and the curve is space-like. Our approach will be based on the theory of Killing vector fields along curves (in the sense of \cite{LS}).

A vector field along a curve $\gamma$, which infinitesimally preserves the arc length parameterization, is said to be a \emph{Killing vector field along $\gamma$} if this curve evolves in the direction of this vector field without changing shape, only position. Let $\gamma:\mathbb{R}\longrightarrow\mathbb{H}_\epsilon^2$ be a convex $p$-elastic curve with nonconstant periodic curvature $\kappa$ and denote by $\{T,N\}$ its associated Frenet frame. Adapting the computations of \cite{LS} (see also \cite{ABG,NEW}), we obtain that the vector field
\begin{equation}\label{J}
	\mathcal{J}:=(p-1)\kappa^p \,T+p\frac{d}{ds}\left(\kappa^{p-1}\right)N\,,
\end{equation}
defined along $\gamma$ is a Killing vector field along the curve.

\begin{rem}
	From the conservation law \eqref{fi}, we obtain that the Killing vector field $\mathcal{J}$ along $\gamma$ is space-like. That is, $\langle\mathcal{J},\mathcal{J}\rangle>0$ for every $s\in\mathbb{R}$.
\end{rem}

By a simple argument involving comparing the corresponding dimensions, it can be shown that Killing vector fields along curves are, precisely, the restrictions to these curves of Killing vector fields defined on the whole space. Therefore, Killing vector fields along curves can uniquely be extended to classical Killing vector fields, which are the infinitesimal generators of isometries. The isometries of $\mathbb{H}_\epsilon^2\subset\mathbb{L}^3$ are those of the Lorentz-Minkowski space $\mathbb{L}^3$ that preserve the quadric $\mathbb{H}_\epsilon^2$. The orbits of these isometries (and so the integral curves of Killing vector fields) may be circles, parabolas or hyperbolas. All these curves have constant curvature.

There are three different types of Killing vector fields in $\mathbb{H}_\epsilon^2$, which depend on the type of their integral curves. We next show that the unique extension of $\mathcal{J}$, \eqref{J}, to $\mathbb{H}_\epsilon^2$ is of \emph{rotational type}, that is, its integral curves are circles. Moreover, we will use this to obtain a explicit parameterization of the $p$-elastic curve $\gamma$.

\begin{theorem}\label{paramthm}
	Let $\gamma:\mathbb{R}\longrightarrow\mathbb{H}_\epsilon^2$ be a convex $p$-elastic curve with nonconstant periodic curvature $\kappa$. Then, the unique extension to $\mathbb{H}_\epsilon^2$ of the Killing vector field $\mathcal{J}$, \eqref{J}, along the $p$-elastic curve is a rotational Killing vector field. Furthermore, the $p$-elastic curve $\gamma$ in $\mathbb{H}_\epsilon^2$ can be parameterized, in terms of its arc length parameter $s\in\mathbb{R}$, as
	\begin{equation}\label{param}
		\gamma(s)=\frac{1}{\sqrt{-a}}\left(\sqrt{\epsilon_2a+p^2\kappa^{2(p-1)}\,}\cos\theta(s),\sqrt{\epsilon_2a+p^2\kappa^{2(p-1)}\,}\sin\theta(s),p\kappa^{p-1}\right),
	\end{equation}
	where
	\begin{equation}\label{angular}
		\theta(s):=\epsilon_2(p-1)\sqrt{-a}\int\frac{\kappa^p}{\epsilon_2a+p^2\kappa^{2(p-1)}}\,ds\,,
	\end{equation}
	is the angular progression and $\epsilon_2=(-1)^\epsilon$.
\end{theorem}
\textit{Proof.} Let $\gamma$ be a convex $p$-elastic curve with nonconstant periodic curvature $\kappa$. Since the curvature $\kappa$ satisfies the Euler-Lagrange equation \eqref{EL}, it is a straightforward computation to check that the arc length parameterization and the curvature are preserved through a variation of $\gamma$ with variational vector field $\mathcal{J}$. This means that $\mathcal{J}$, \eqref{J}, is a Killing vector field along $\gamma$ (for details, see \cite{ABG,LS} and references therein).

Then, $\mathcal{J}$ can be uniquely extended to a vector field in $\mathbb{H}_\epsilon^2$, that is, denoting this extension also by $\mathcal{J}$, $\mathcal{J}\in\mathfrak{X}(\mathbb{H}_\epsilon^2)$. Moreover, $\mathcal{J}\in\mathfrak{X}(\mathbb{H}_\epsilon^2)$ is a Killing vector field in the classical sense. We will next show that the integral curves of $\mathcal{J}\in\mathfrak{X}(\mathbb{H}_\epsilon^2)$ are circles and, hence, $\mathcal{J}$ is rotational.

For this purpose, we notice that since $\kappa$ is periodic, it must attain its maximum. Say that $s_o\in\mathbb{R}$ is one value of the arc length parameter such that $\kappa(s_o)=\alpha$ is the maximum of $\kappa$. It then follows that $\kappa'(s_o)=0$, i.e., the point $p_o=\gamma(s_o)$ is a \emph{vertex}. We now compute the curvature of the integral curve $\widetilde{\gamma}$ of $\mathcal{J}$ passing through $p_o=\gamma(s_o)$. Since $\widetilde{\gamma}$ is an integral curve of $\mathcal{J}$, its unit tangent is 
$$\widetilde{T}=\frac{\mathcal{J}}{\lVert\mathcal{J}\rVert}\,,$$
where $\lVert \mathcal{J}\rVert=\sqrt{\langle\mathcal{J},\mathcal{J}\rangle}$ (recall that $\mathcal{J}$ is space-like and so the square root of $\langle\mathcal{J},\mathcal{J}\rangle$ is well defined). Differentiating $\mathcal{J}$, \eqref{J}, with respect to the arc length parameter $s\in\mathbb{R}$ and evaluating at $s_o$, we have
$$\nabla_T\mathcal{J}(s_o)=\left(\epsilon_2(p-1)\alpha^{p+1}+p\frac{d^2}{ds^2}_{\lvert_{s=s_o}}\left(\kappa^{p-1}\right)\right)N(s_o)=\epsilon_2\,p\,\alpha^{p-1}N(s_o)\,,$$
since $\kappa(s_o)=\alpha$ and $\kappa'(s_o)=0$. Observe that we have used the Euler-Lagrange equation \eqref{EL} to obtain $\kappa''(s_o)$. From \eqref{J}, it follows that at $s_o$, $\mathcal{J}(s_o)=(p-1)\alpha^pT(s_o)$ and so the unit tangent to $\widetilde{\gamma}$ there is $\widetilde{T}=T$. Hence, at $s_o$, we have
$$\nabla_{\widetilde{T}}\widetilde{T}=\nabla_T\left(\frac{\mathcal{J}}{\lVert\mathcal{J}\rVert}\right)=\frac{1}{\lVert\mathcal{J}\rVert}\nabla_T\mathcal{J}=\frac{\epsilon_2p}{(p-1)\alpha}N\,,$$
since $T\langle\mathcal{J},\mathcal{J}\rangle=2\langle\nabla_T\mathcal{J},\mathcal{J}\rangle=0$ holds. We then deduce that the curvature of $\widetilde{\gamma}$ at $s_o$ is
$$\widetilde{\kappa}(s_o)=\frac{p}{(p-1)\alpha}\,.$$
Moreover, since $\widetilde{\gamma}$ is an integral curve of the Killing vector field $\mathcal{J}\in\mathfrak{X}(\mathbb{H}_\epsilon^2)$, it has constant curvature $\widetilde{\kappa}_o=\widetilde{\kappa}(s_o)$. Applying the Gauss Formula and the conservation law \eqref{fi} evaluated at $s_o$, we obtain
$$\langle\overline{\nabla}_{\widetilde{T}}\widetilde{T},\overline{\nabla}_{\widetilde{T}}\widetilde{T}\rangle=\frac{\epsilon_2p^2}{(p-1)^2\alpha^2}-\epsilon_2=-\frac{a}{(p-1)^2\alpha^{2p}}>0\,,$$
because $a<0$ must be satisfied in order the $p$-elastic curve $\gamma$ to have periodic curvature $\kappa$. This shows that the integral curve of $\mathcal{J}$ passing through $p_o=\gamma(s_o)$ is a circle of $\mathbb{L}^3$ embedded in $\mathbb{H}_\epsilon^2$. Therefore, all the integral curves of $\mathcal{J}$ are circles and, hence, it is a rotational Killing vector field of $\mathbb{H}_\epsilon^2$.

The fact that the vector field $\mathcal{J}$ is rotational suggests that the parameterization of $\gamma$ must involve trigonometric functions. Indeed, we will show that \eqref{param} is an arc length parameterization of the $p$-elastic curve $\gamma$.

First observe that from \eqref{param}, $\langle\gamma(s),\gamma(s)\rangle=-\epsilon_2$ and, hence, $\gamma(\mathbb{R})\subset\mathbb{H}_\epsilon^2$. Moreover, differentiating \eqref{param} and using the definition of the angular progression $\theta(s)$, we obtain $\langle\gamma'(s),\gamma'(s)\rangle=1$. This means that \eqref{param} is an arc length parameterized space-like curve. Finally, differentiating \eqref{param} twice and employing the Gauss Formula, we deduce that the curvature of the parameterized curve \eqref{param} is, precisely, $\kappa(s)$. Therefore, \eqref{param} is an arc length parameterization of the $p$-elastic curve $\gamma$, up to rigid motions of $\mathbb{H}_\epsilon^2$. \hfill$\square$

\begin{rem}
	The explicit parameterizations \eqref{param} of $p$-elastic curves have been found adapting the computations of \cite{ABG}. This approach involves employing the rotational Killing vector field in an essential way. The parameterizations \eqref{param} may as well be found with an approach based on the Marsden-Weinstein reduction method.
\end{rem}

Some geometric properties can be directly deduced from the parameterization \eqref{param}. We summarize them in the following result.

\begin{proposition}\label{properties}
	Let $\gamma:\mathbb{R}\longrightarrow\mathbb{H}_\epsilon^2$ be a convex $p$-elastic curve with nonconstant periodic curvature $\kappa$. Then, up to rigid motions of $\mathbb{H}_\epsilon^2$:
	\begin{enumerate}[(i)]
		\item The trajectory of $\gamma$ is contained between two parallels of $\mathbb{H}_\epsilon^2$. If $\gamma\subset\mathbb{H}_0^2$ is a hyperbolic curve, it never meets the pole $(0,0,1)$. If $\gamma\subset\mathbb{H}_1^2$ is a pseudo-hyperbolic curve, then it is entirely contained in $\mathbb{H}_1^2\lvert_-=\mathbb{H}_1^2\cap\{z<0\}$
		\item The curve $\gamma$ meets the bounding parallels tangentially at the maximum and minimum values of its curvature, $\alpha\equiv\alpha(a)$ and $\beta\equiv\beta(a)$, respectively. Indeed, these parallels, located at heights $z=p\alpha^{p-1}/\sqrt{-a}$ and $z=p\beta^{p-1}/\sqrt{-a}$ respectively, are the integral curves of the Killing vector field $\mathcal{J}$ passing through the corresponding points.
		\item The angular progression \eqref{angular} is monotonic with respect to the arc length parameter of the curve. This means that the trajectory of $\gamma$ winds around the $z$-axis, without moving backwards.
		\item The $p$-elastic curve $\gamma$ is closed if and only if the angular progression \eqref{angular} along a period $\varrho\equiv\varrho(a)$ of the curvature is a rational multiple of $2\pi$. That is, if and only if,
		\begin{equation}\label{Lambda}
			\Lambda_p(a):=\epsilon_2(p-1)\sqrt{-a}\int_0^\varrho \frac{\kappa^p}{\epsilon_2a+p^2\kappa^{2(p-1)}}\,ds=2\pi\,q\,,
		\end{equation}
		for some $q\in\mathbb{Q}$.
	\end{enumerate}
\end{proposition}
\textit{Proof.} Assume that $\gamma$ is a convex $p$-elastic curve in $\mathbb{H}_\epsilon^2$ with nonconstant periodic curvature $\kappa$ parameterized by \eqref{param}. Since the curvature is periodic it attains both its maximum and minimum values, namely, $\alpha$ and $\beta$. At the points where those values are attained $\gamma$ is tangent to the integral curves of $\mathcal{J}$ (see the proof of Theorem \ref{paramthm}), which are parallel circles. These parallels are located at heights $z=p\alpha^{p-1}/\sqrt{-a}$ and $z=p\beta^{p-1}/\sqrt{-a}$, respectively. It follows from $\beta\leq \kappa\leq \alpha$ that $\gamma$ is contained in the region bounded by those parallels. In fact, if $\gamma$ is a hyperbolic curve $(\epsilon=0)$ then $p>1$, while if $\gamma$ is a pseudo-hyperbolic curve $(\epsilon=1)$ then $p<0$. In any case,
$$\frac{1}{\sqrt{-a}}\,p\,\beta^{p-1}\leq \frac{1}{\sqrt{-a}}\,p\,\kappa^{p-1}=z(s)\leq \frac{1}{\sqrt{-a}}\,p\,\alpha^{p-1}\,.$$
This suffices to show that in the case that $\gamma$ is a pseudo-hyperbolic curve $(\epsilon=1)$, then the curve satisfies $z<0$ everywhere, since $p<0$ and $\alpha>\beta>0$ holds. In the case that $\gamma$ is a hyperbolic curve $(\epsilon=0)$, none of the above parallels degenerate to the pole $(0,0,1)$ so the curve $\gamma$ cannot pass through this point. Assume to the contrary that there exist $s_o\in\mathbb{R}$ such that $\gamma(s_o)=(0,0,1)$. Then, from \eqref{param} we deduce that $a+p^2\kappa^{2(p-1)}=0$. However, notice that from the conservation law \eqref{fi},
$$\langle\mathcal{J},\mathcal{J}\rangle=(p-1)^2\kappa^{2p}+p^2\left(\frac{d}{ds}\left(\kappa^{p-1}\right)\right)^2=a+p^2\kappa^{2(p-1)}\,,$$
contradicting the fact that $\mathcal{J}$ is space-like (which, by the way, is a consequence of the middle expression being positive). This shows the first and second properties.

For the third property, we just differentiate the angular progression $\theta(s)$, \eqref{angular}, with respect to the arc length parameter $s$, obtaining
$$\theta'(s)=\epsilon_2(p-1)\sqrt{-a}\,\frac{\kappa^p}{\epsilon_2a+p^2\kappa^{2(p-1)}}\neq 0\,,$$
because $\gamma$ is convex.

Finally, it is clear that, since the curves $\gamma$ wind around the $z$-axis, they will close if and only if $\theta(\varrho)-\theta(0)$ is a rational multiple of $2\pi$. Here, $\varrho$ denotes the least period of the curvature $\kappa$. The condition $\theta(\varrho)-\theta(0)$ can be rewritten as \eqref{Lambda}. \hfill$\square$
\\

\begin{figure}[h!]
	\centering
	\includegraphics[height=5.55cm]{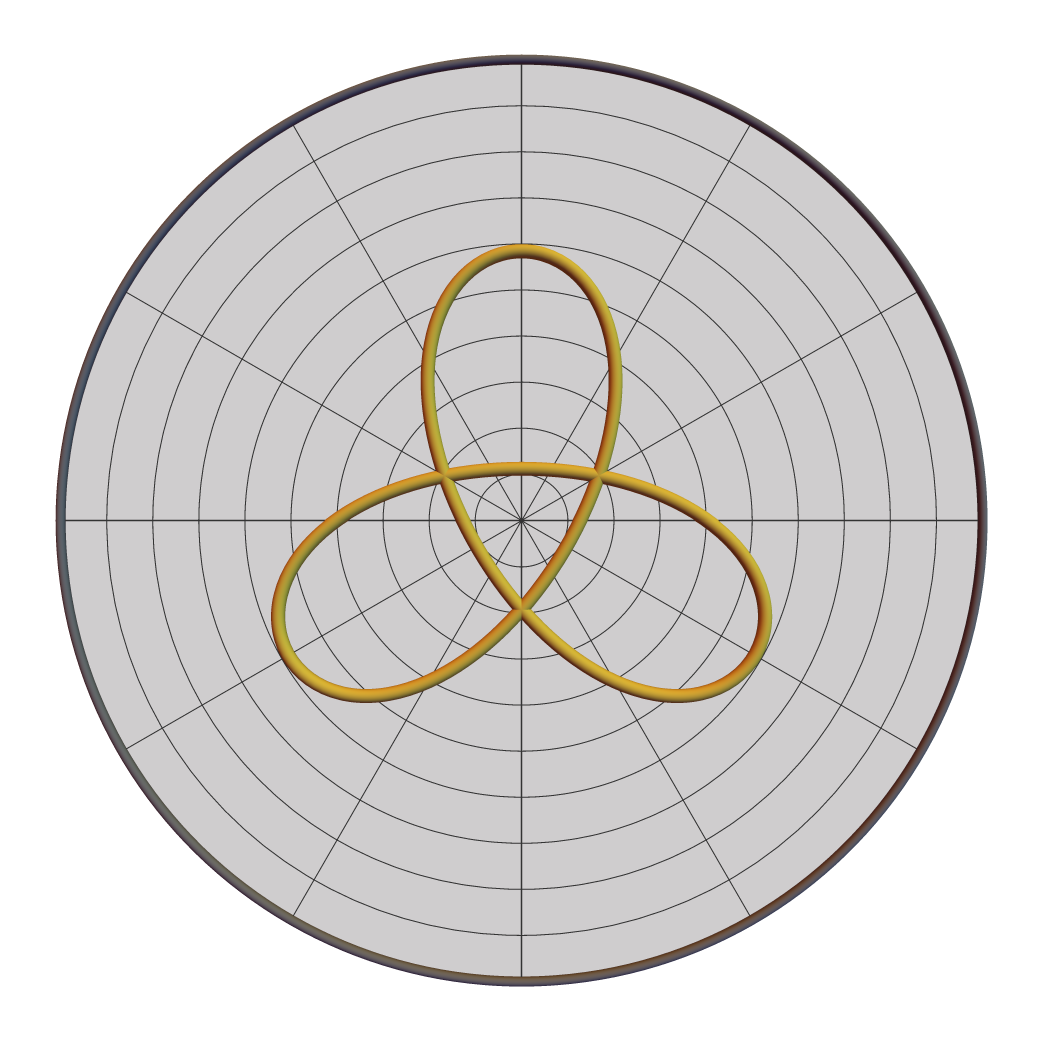}
	\includegraphics[height=5.55cm]{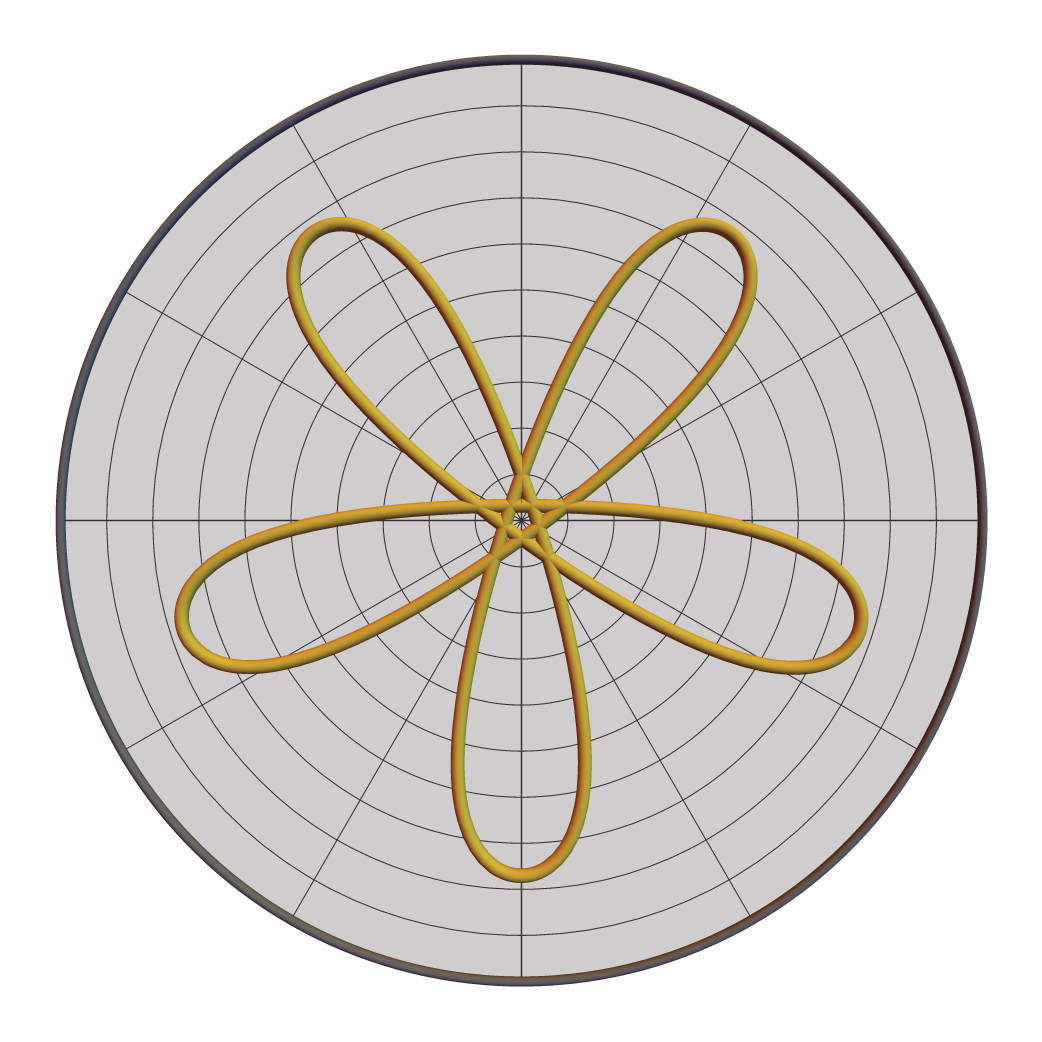}
	\includegraphics[height=5.55cm]{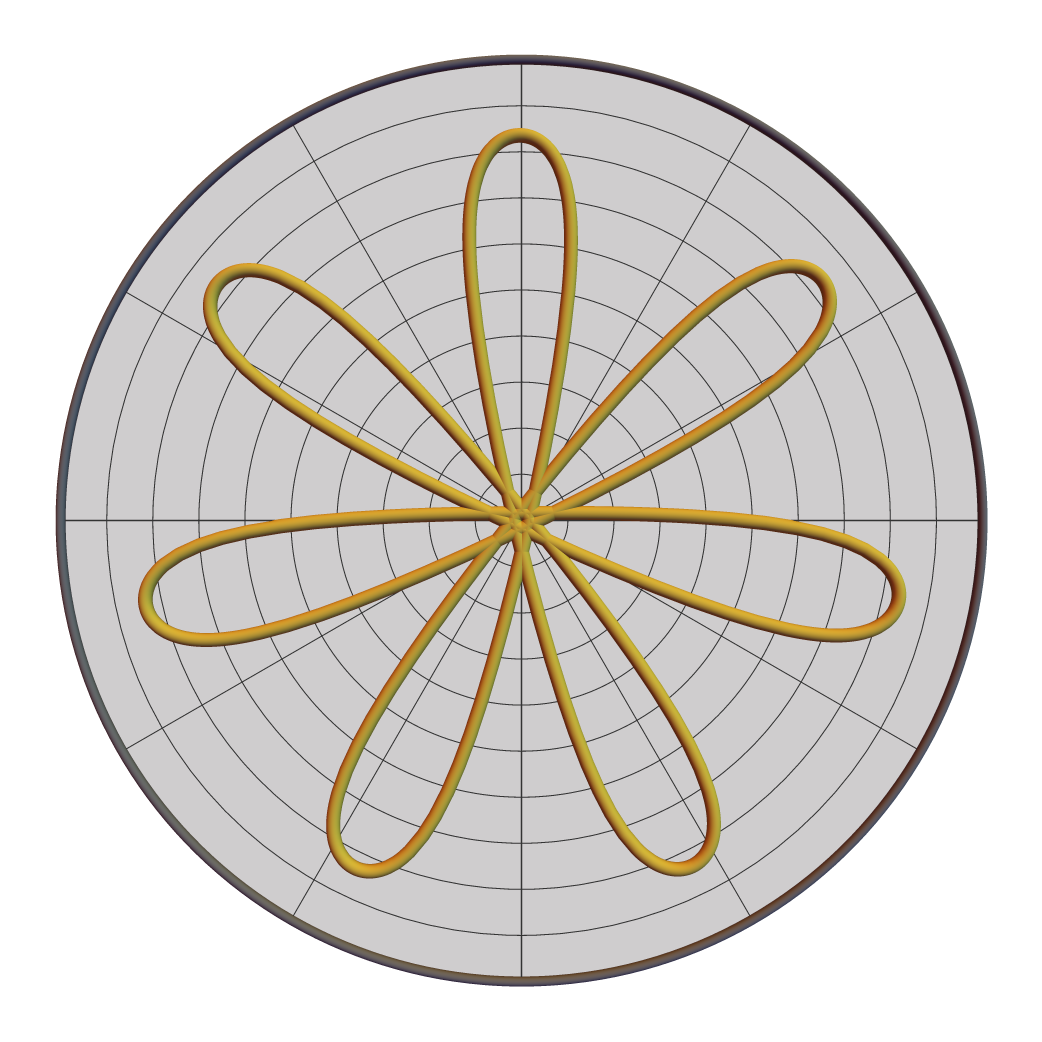}
	\caption{Three hyperbolic $p$-elastic curves for $p=3/2$ in $\mathbb{H}^2_0$ corresponding to the values $q=2/3$, $q=3/5$ and $q=4/7$, respectively. They are represented in the Poincar\'e disk model.}
	\label{F1}
\end{figure}

\begin{figure}[h!]
	\centering
	\includegraphics[height=5.55cm]{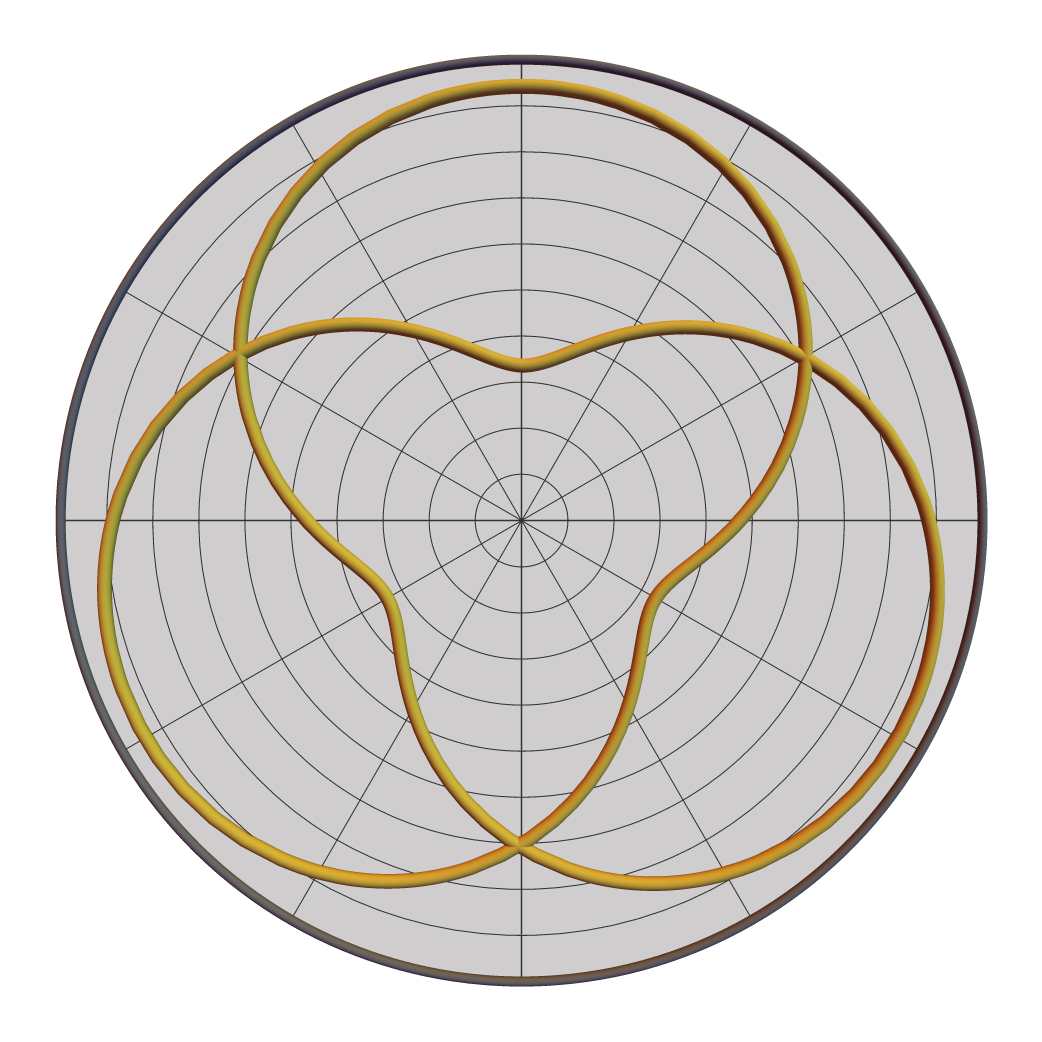}
	\includegraphics[height=5.55cm]{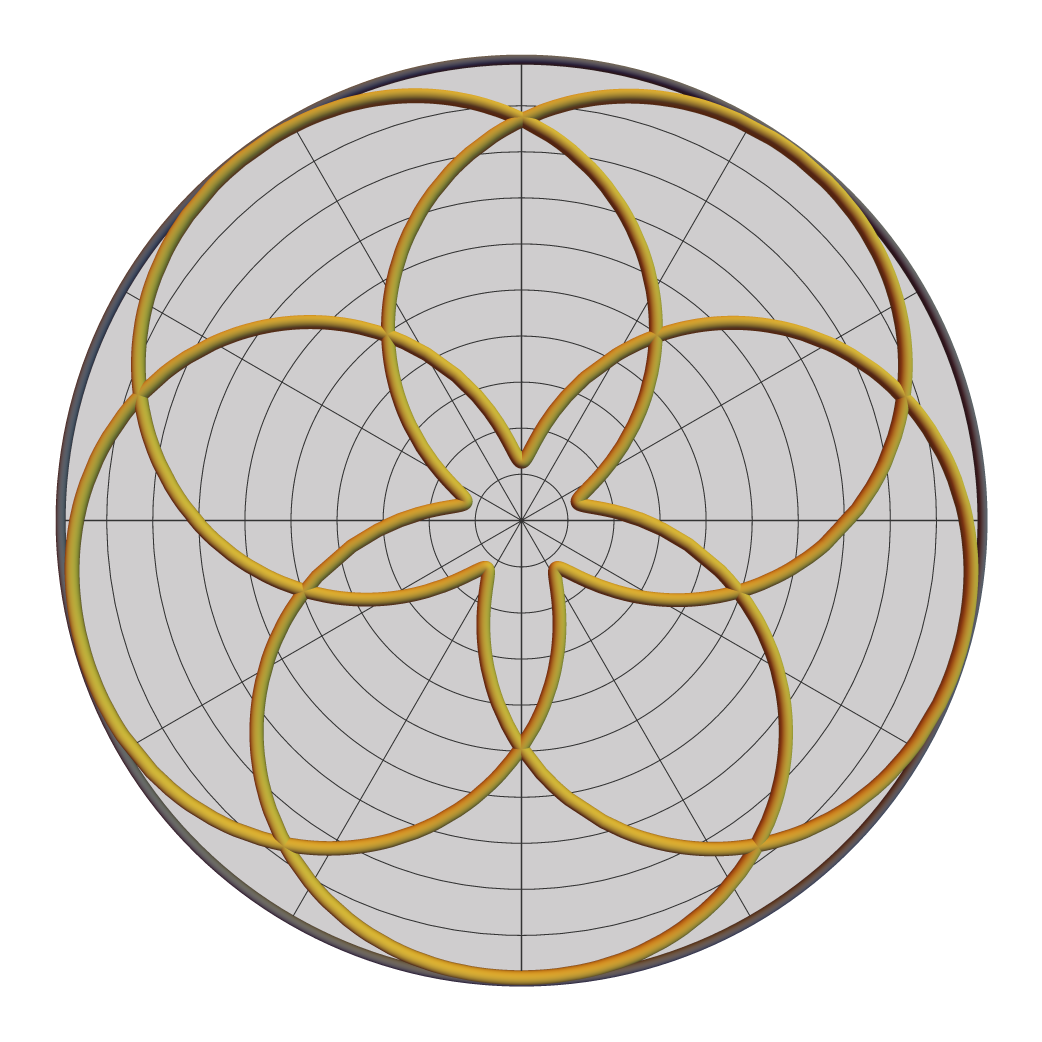}
	\includegraphics[height=5.55cm]{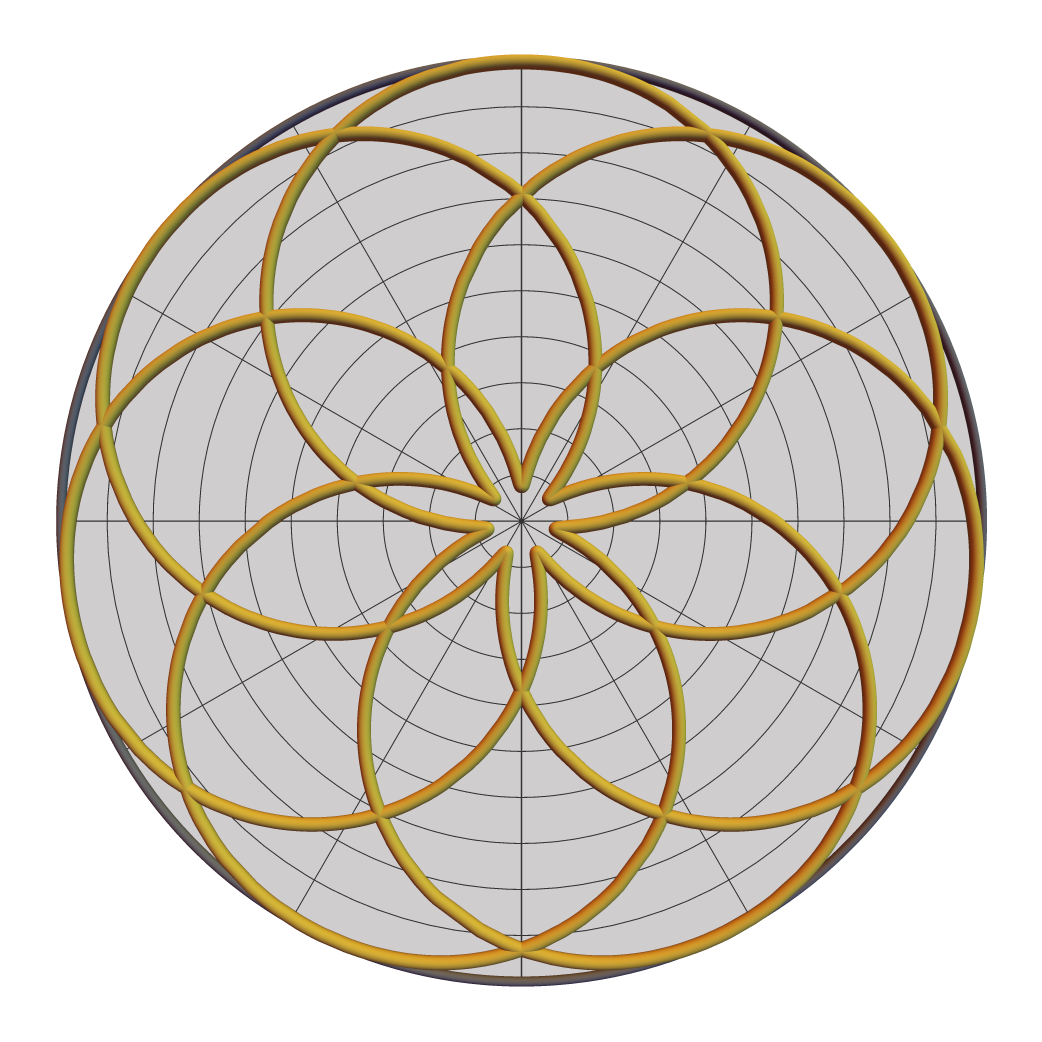}
	\caption{Three (space-like) pseudo-hyperbolic $p$-elastic curves for $p=-1$ in $\mathbb{H}_1^2$ corresponding to the values $q=2/3$, $q=3/5$ and $q=4/7$, respectively. They are represented in the once-punctured unit disk $\mathring{\mathbb{D}}$. (Observe that the transformation that takes $\mathbb{H}_1^2\lvert_-$ to $\mathring{\mathbb{D}}$ involves an inversion.)}
	\label{F2}
\end{figure}

In the following figures we illustrate the geometric properties described in Proposition \ref{properties}. Figure \ref{F1} shows three (hyperbolic) $p$-elastic curves in $\mathbb{H}_0^2$ for the value $p=3/2$, while in Figure \ref{F2} we represent three (pseudo-hyperbolic) $p$-elastic curves in $\mathbb{H}_1^2$ for the value $p=-1$. These closed curves $\gamma$ correspond with the value $0>a>a_*$ such that $\Lambda_p(a)$ is the rational multiple of $2\pi$ given by $q=2/3$, $q=3/5$ and $q=4/7$, respectively. The curves are shown in the unit disks $\mathbb{D}$ and $\mathring{\mathbb{D}}$, respectively, as explained in Section 2.

\section{Periodic Critical Curves}

In Proposition \ref{properties} we obtained a sufficient condition for a convex $p$-elastic curve with nonconstant periodic curvature to be closed. This sufficient condition requires the analysis of $\Lambda_p(a)$, \eqref{Lambda}, as a function of the constant of integration $0>a>a_*$. In this section we will study the existence of closed $p$-elastic curves, for suitable values of $p\in\mathbb{R}$, immersed in both $\mathbb{H}_0^2$ and $\mathbb{H}_1^2$. In addition, we will see that closed $p$-elastic curves are encoded by two natural parameters under suitable restrictions. Moreover, we will obtain a simple expression to compute their energy and discuss about their evolution with respect to the energy parameter $p$.

\subsection{Existence}

Let $\gamma:\mathbb{R}\longrightarrow\mathbb{H}_\epsilon^2$ be a convex $p$-elastic curve in $\mathbb{H}_\epsilon^2$ with nonconstant periodic curvature $\kappa$. We recall that the periodicity of $\kappa$ imposes some restrictions on the parameter $p\in\mathbb{R}$, in addition to $0>a>a_*$ and $\gamma$ being space-like (Theorem \ref{periodic}). In fact, if $\gamma$ is a hyperbolic curve $(\epsilon=0)$ then $p>1$ must hold, and if $\gamma$ is a pseudo-hyperbolic curve $(\epsilon=1)$ then $p<0$ must hold. It then follows from Proposition \ref{properties} that the $p$-elastic curve $\gamma$ is closed if and only if the function $\Lambda_p:(a_*,0)\subset\mathbb{R}\longrightarrow \mathbb{R}$ defined in \eqref{Lambda} satisfies $\Lambda_p(a)=2\pi q$ for some $q\in\mathbb{Q}$. In order to check this, it is convenient to use the conservation law \eqref{fi} to make a change of variable so that $\Lambda_p$ can be rewritten as
\begin{equation}\label{terrible}
	\Lambda_p(a)=2p(p-1)^2\sqrt{-a}\int_\beta^\alpha \frac{\kappa^{2(p-1)}}{(a+\epsilon_2p^2\kappa^{2(p-1)})\sqrt{a-\epsilon_2(p-1)^2\kappa^{2p}+\epsilon_2p^2\kappa^{2(p-1)}\,}}\,d\kappa\,,
\end{equation}
where $\epsilon_2=(-1)^\epsilon$ and $0<\beta<\alpha$ are the only two positive solutions of
\begin{equation}\label{solutions}
	f_{p,a}(\kappa):=a-\epsilon_2(p-1)^2\kappa^{2p}+\epsilon_2p^2\kappa^{2(p-1)}=0\,.
\end{equation}
The assertion about the positive solutions $\beta<\alpha$ follows from the analysis of the functions $Q_{p,a}$ and $\widetilde{Q}_{p,a}$ carried out in the proof of Theorem \ref{periodic}. Indeed, the function $f_{p,a}$ given in \eqref{solutions} is, precisely, $\widetilde{Q}_{p,a}$ and, hence, a suitable multiple of $Q_{p,a}$ too. 

In the following lemma we show that the function $\Lambda_p$ is continuous at every $a\in(a_*,0)$ and that $\Lambda_p(a)\to\sqrt{2}\,\pi$ as $a\to a_*^+$.

\begin{lemma}\label{lemma}
	Assume that $\epsilon=0$ holds when $p>1$ and that $\epsilon=1$ holds when $p<0$. Then, for every $p\in(-\infty,0)\cup(1,\infty)$ the function $\Lambda_p:(a_*,0)\subset\mathbb{R}\longrightarrow\mathbb{R}$ given by \eqref{terrible} is continuous and
	$$\lim_{a\to a_*^+} \Lambda_p(a)=\sqrt{2}\,\pi\,,$$
	holds.
\end{lemma}
\textit{Proof.} Fix an admissible value of $p\in\mathbb{R}$ and consider the corresponding value of $\epsilon=0,1$ such that the restrictions of the statement are satisfied. From Theorem \ref{periodic} we conclude that when $\kappa\in(\beta,\alpha)$, $f_{p,a}(\kappa)$ is positive and, hence, the square root arising in \eqref{terrible} is well defined. Furthermore, in the case $\epsilon=0$, the same reasoning shows that
$$a+\epsilon_2p^2\kappa^{2(p-1)}>\epsilon_2(p-1)^2\kappa^{2p}>0\,,$$
since $\epsilon_2=(-1)^\epsilon=1$. In the case $\epsilon=1$, $\epsilon_2=(-1)^\epsilon=-1$ and, hence, it is clear that $a+\epsilon_2p^2\kappa^{2(p-1)}<0$. In both cases, the integrand of \eqref{terrible} is continuous and so is the function $\Lambda_p$.

At this point, we will compute the limit of $\Lambda_p$ when $a\to a_*^+$. For this purpose, we will employ the result of Corollary 4.2 of \cite{Pe}. Recall that the function $f_{p,a}$ defined in \eqref{solutions} has a local maximum at $\kappa=\kappa_*$ (see the proof of Theorem 4.1). Thus, the conditions of Corollary 4.2 of \cite{Pe} are satisfied and
$$\lim_{a\to a_*^+} \Lambda_p(a)=2p(p-1)^2\sqrt{-a_*}\frac{\kappa_*^{2(p-1)}}{\left(a_*+\epsilon_2p^2\kappa_*^{2(p-1)}\right)\sqrt{2\epsilon_2p^2(p-1)\kappa_*^{2(p-2)}}}\,\pi=\sqrt{2}\,\pi\,.$$
This finishes the proof. \hfill$\square$
\\

Since $\Lambda_p$ is continuous, to check that there exist rational numbers $q\in\mathbb{Q}$ such that $\Lambda_p(a)=2\pi q$, it suffices to show that the image of $\Lambda_p$ contains an interval of the real line. For the hyperbolic plane $\mathbb{H}_0^2$ (since $\epsilon=0$, then $p>1$ must hold) we will next prove that the image of $\Lambda_p$ contains the interval $(\pi,\sqrt{2}\,\pi)$ and, hence, for every $q\in\mathbb{Q}$ such that $1<2q<\sqrt{2}$ there exists a closed $p$-elastic curve in $\mathbb{H}_0^2$ with nonconstant curvature.

\begin{theorem}\label{closed}
	For every $p>1$ and rational number $q\in\mathbb{Q}$ satisfying $1<2q<\sqrt{2}$, there exists a closed $p$-elastic curve in the hyperbolic plane $\mathbb{H}_0^2$ with nonconstant curvature.
\end{theorem}
\textit{Proof.} Let $p>1$ and $\epsilon=0$ be fixed. From Lemma \ref{lemma} we know that $\Lambda_p$ is a continuous function satisfying $\Lambda_p\to \sqrt{2}\,\pi$ as $a\to a_*^+$. By computing the limit of $\Lambda_p$ as $a\to 0^-$ we will show that the image of $\Lambda_p$ contains the interval $(\pi,\sqrt{2}\,\pi)$.

To compute the limit of $\Lambda_p$ when $a\to 0^-$, we will proceed as in \cite{MOP} (see also \cite{GPT} and references therein) and argue by means of the Cauchy's Integral Theorem. We will compute this limit for a rational number $p=l/t\in\mathbb{Q}$ where $l,t\in\mathbb{N}$ and $t<l$ are two positive relatively prime integers. The result for general values of $p>1$ will then follow from the continuity of $\Lambda_p$ in $p$.

Assume that $p=l/t\in\mathbb{Q}$ and consider the change of variable $\kappa^2=u^t$ to rewrite \eqref{terrible} as
$$\Lambda_{p=l/t}(a)=l(p-1)^2\sqrt{-a}\int_{\widehat{\beta}}^{\widehat{\alpha}} \frac{u^{l-t/2-1}}{\left(a+p^2u^{l-t}\right)\sqrt{a-(p-1)^2u^l+p^2u^{l-t}}}\,du\,,$$
where $\widehat{\alpha}=\alpha^{2/t}$ and $\widehat{\beta}=\beta^{2/t}$. The advantage of this new expression is that now what is inside the square root, that is $f_{p,a}(u)$, is a polynomial of degree $l$. Moreover, this polynomial has exactly two positive roots, namely, $\widehat{\alpha}>\widehat{\beta}$ and, in addition, it follows from Descartes' rule of signs that it has, at most, one negative root. This negative root, which we denote by $\delta$, exists when $l$ is odd. The rest of the roots lie in $\mathbb{C}\setminus\mathbb{R}$ and we denote them by $\omega_j$ and $\overline{\omega}_j$, $j=1,...,\lfloor (l-2)/2 \rfloor$, where the upper line denotes the complex conjugate. The polynomial $a+p^2 u^{l-t}$ has just one positive root, say $u_o$, which satisfies $u_o<\widehat{\beta}$. The other roots are complex, namely, $u_i$ and $\overline{u}_i$, $i=1,...,\lfloor (l-t-1)/2\rfloor$; and, when $l-t$ is even there also exists one negative root which is, precisely, $-u_o$.

We then define the complex function $h(z)$ by
$$h(z)=l(p-1)^2\sqrt{-a}\frac{z^{l-1}}{(a+p^2z^{l-t})\sqrt{z(a-(p-1)^2z^l+p^2z^{l-t})}}\,,$$
where the complex square root is defined as $z=\sqrt{r\,e^{i\theta}}=\sqrt{\lvert r\rvert}\,e^{i\theta/2}$. The function $h(z)$ defined above is well defined and analytic far from the non positive part of the $u$-axis, the roots of the polynomial $f_{p,a}(u)$, and the roots of $a+p^2u^{l-t}$.

Let $\sigma$ be a big circle enclosing all the singularities of $h$ (with a small cut so that the non positive part of the $u$-axis is not inside $\sigma$), $\sigma_*$ a small path around $\widehat{\alpha}$ and $\widehat{\beta}$, and $\sigma_\circ$ a little circle around each of the rest of the singularities of $h$ that lie in $\mathbb{C}\setminus\mathbb{R}^-$. All the paths are assumed to be oriented counter-clockwise. Observe that along $\sigma_*$ we have, from the definition of the complex square root,
\begin{eqnarray*}
	\lim_{\epsilon\to 0^+} h(u+i\epsilon)&=&h(u)\,,\\
	\lim_{\epsilon\to 0^-} h(u+i\epsilon)&=&-h(u)\,,
\end{eqnarray*} 
for every $u\in(\widehat{\beta},\widehat{\alpha})$, since when $\epsilon\to 0^+$ then $\theta\to 0$ in the definition of the square root, while when $\epsilon\to 0^-$, then $\theta\to 2\pi$.

From the analyticity of the function $h$, it follows that
\begin{eqnarray*}
	\int_\sigma h(z)\,dz&=&\int_{\sigma_*} h(z)\,dz+\int_{\sigma_{u_o}} h(z)\,dz\\
	&&+\sum_{i=1}^{\lfloor (l-t-1)/2\rfloor} \left(\int_{\sigma_{u_i}}h(z)\,dz+\int_{\sigma_{\overline{u}_i}}h(z)\,dz\right)\\
	&&+\sum_{j=1}^{\lfloor (l-2)/2\rfloor} \left(\int_{\sigma_{\omega_j}} h(z)\,dz+\int_{\sigma_{\overline{\omega}_j}} h(z)\,dz\right).
\end{eqnarray*} 
We next deduce from Cauchy's Integral Formula that the integrals on the last line above are all zero, while those in the second line have opposite signs for each $i=1,...,\lfloor (l-t-1)/2 \rfloor$, and so they cancel out. In the same way, we compute
$$\int_{\sigma_{u_o}} h(z)\,dz=2\pi i \left(\frac{1}{2\pi i}\int_{\sigma_{u_o}} \frac{g(z)}{z-u_o}\,dz\right)=2\pi i g(u_o)=2\pi\,,$$
where $g(z)=(z-u_o)h(z)$. The last equality follows after a straightforward simplification of $g(u_o)$. Therefore, we conclude with
$$\int_{\sigma_*}h(z)\,dz=\int_\sigma h(z)\,dz -2\pi\,,$$
for every $a\in(a_*,0)$.

Finally, we notice that along the path $\sigma$, $h(z)\to 0$ as $a\to 0^-$. So, using the above limits for $h(z)$ as $\epsilon\to 0$, we obtain
$$\lim_{a\to 0^-} \Lambda_{l/t}(a)=-\frac{1}{2} \lim_{\epsilon\to 0}\left(\lim_{a\to 0^-}\int_{\sigma_*}h(z)\,dz\right)=-\frac{1}{2}\left(-2\pi\right)=\pi\,.$$

Since the image of $\Lambda_p$ contains the interval $(\pi,\sqrt{2}\,\pi)$, it follows that for every rational number $q\in\mathbb{Q}$ such that $2\pi q\in(\pi,\sqrt{2}\,\pi)$ there exists a value $a_q\in(a_*,0)$ satisfying $\Lambda_p(a_q)=2\pi q$. We then deduce from Proposition \ref{properties} that the curve $\gamma\equiv\gamma_{a_q}$ is a closed $p$-elastic curve for any fixed $p>1$. Observe that $2\pi q\in(\pi,\sqrt{2}\,\pi)$ is equivalent to $1<2q<\sqrt{2}$. \hfill$\square$
\\

We next rewrite Theorem \ref{closed} in terms of two natural parameters. Let $\gamma_q$ be a closed $p$-elastic curve in $\mathbb{H}_0^2$ with nonconstant curvature associated to $q\in\mathbb{Q}$. Write the rational number $q\in\mathbb{Q}$ as $q=n/m$ for relatively prime natural numbers $n$ and $m$. Then, $n$ is the winding number of the curve $\gamma_q\equiv\gamma_{n,m}$ and $\gamma_{n,m}$ has an $m$-fold symmetry. In the disk model, $n$ represents the number of times the curve goes around the center of the disk, while $m$ represents the number of lobes of the curve.

\begin{corollary}\label{cor}
	For every $p>1$ and pair of relatively prime natural numbers $(n,m)$ satisfying $m<2n<\sqrt{2}\,m$, there exists a closed $p$-elastic curve in the hyperbolic plane $\mathbb{H}_0^2$ with nonconstant curvature.
\end{corollary}

\begin{figure}[h!]
	\centering
	\includegraphics[height=5.55cm]{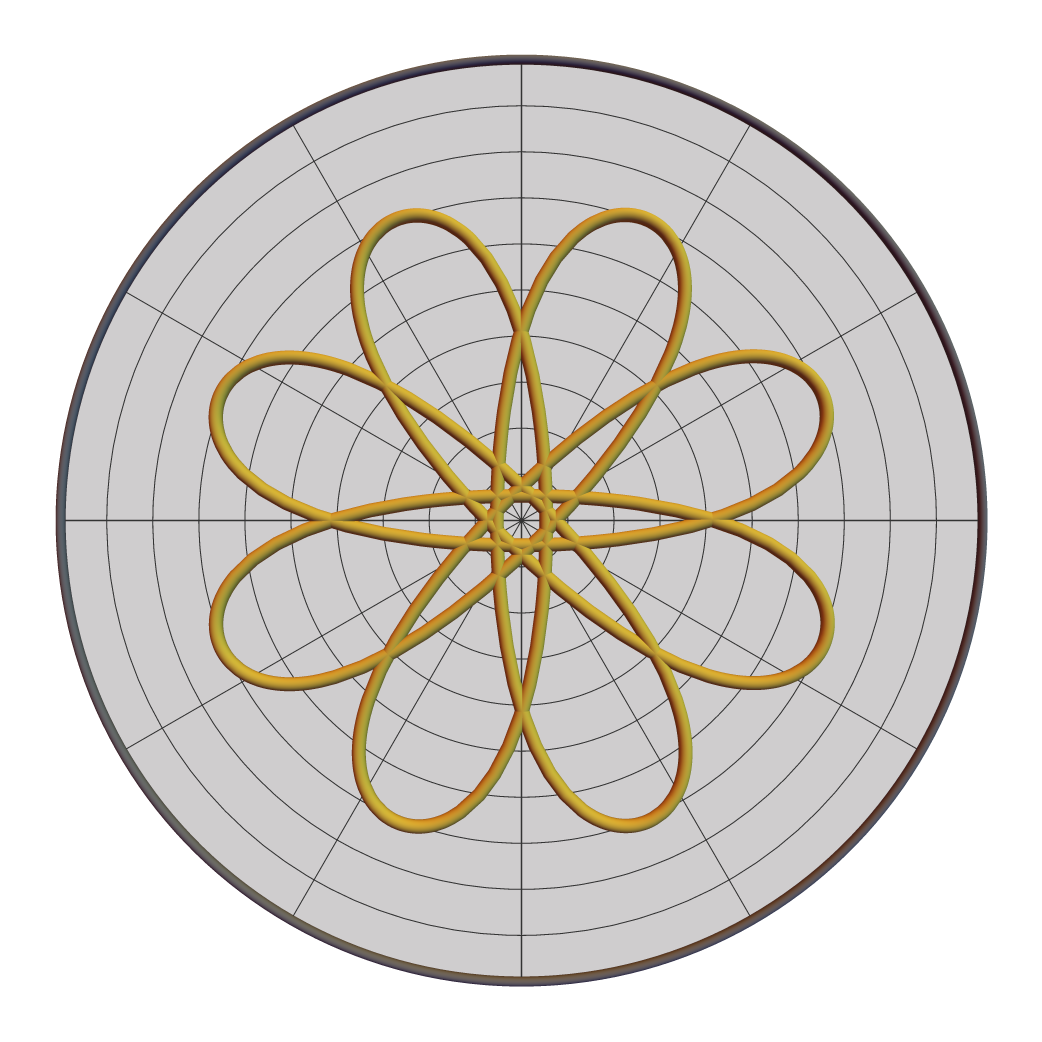}
	\includegraphics[height=5.55cm]{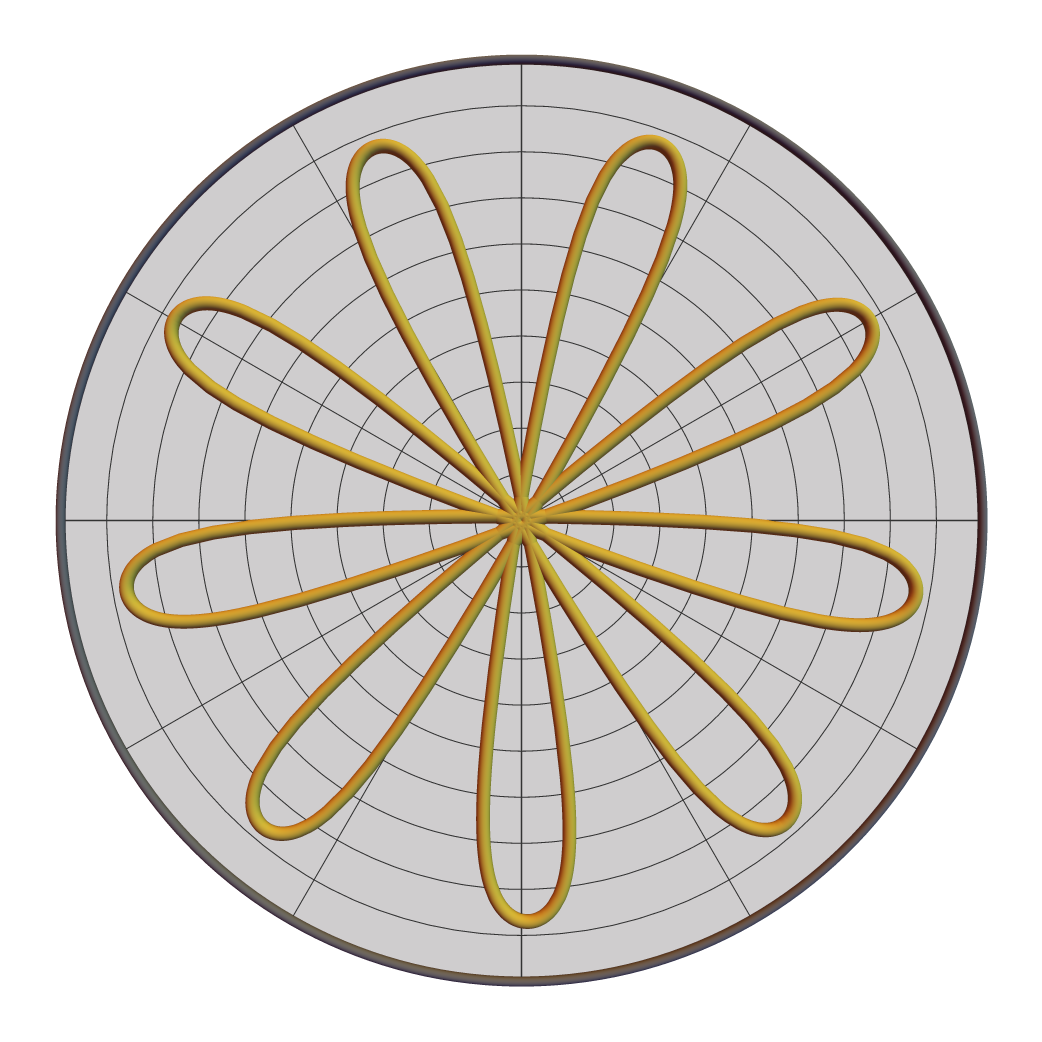}
	\includegraphics[height=5.55cm]{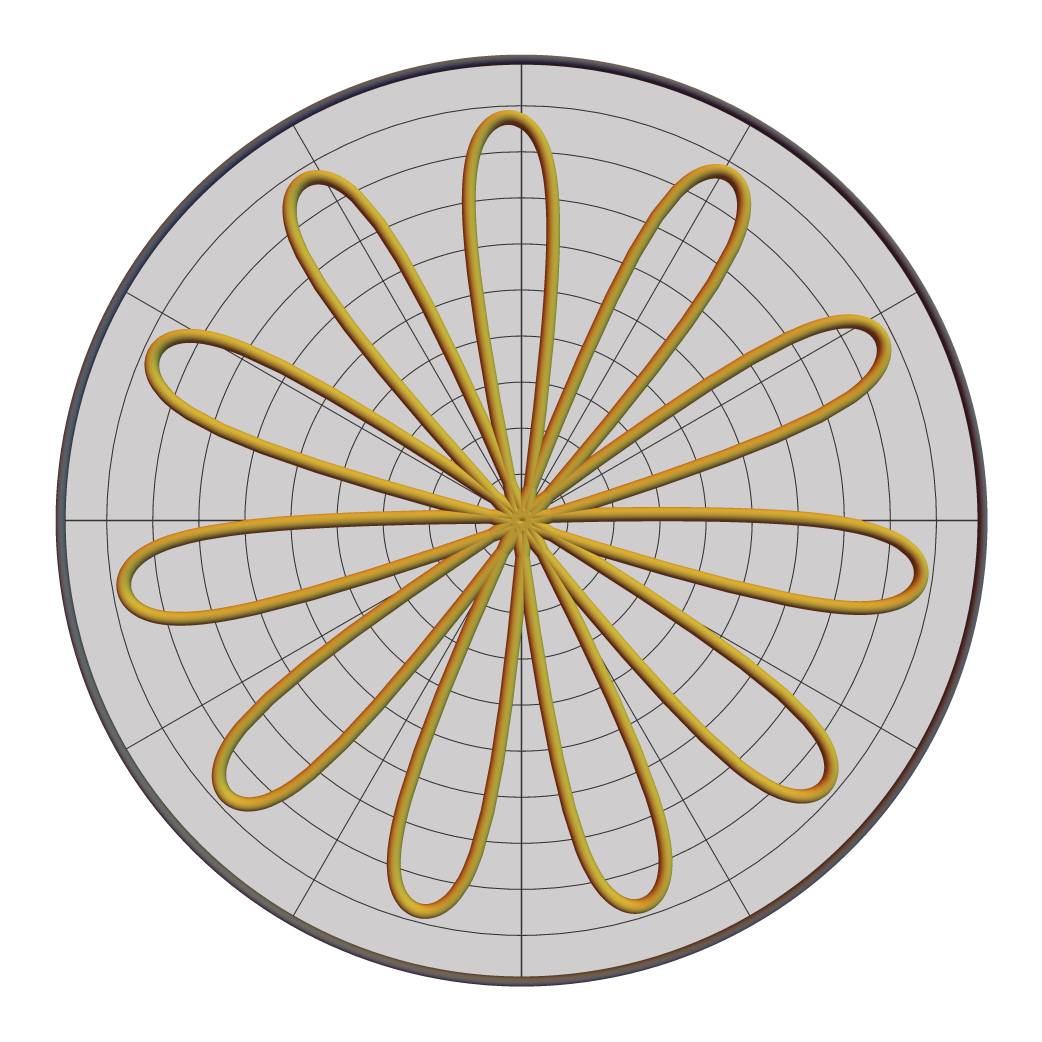}
	\caption{Closed hyperbolic $p$-elastic curves for $p=3/2$ in $\mathbb{H}^2_0$ corresponding to the values $q=5/8$, $q=5/9$ and $q=6/11$, respectively. They are represented in the Poincar\'e disk model.}
	\label{F3}
\end{figure}

\begin{figure}[h!]
	\centering
	\includegraphics[height=5.55cm]{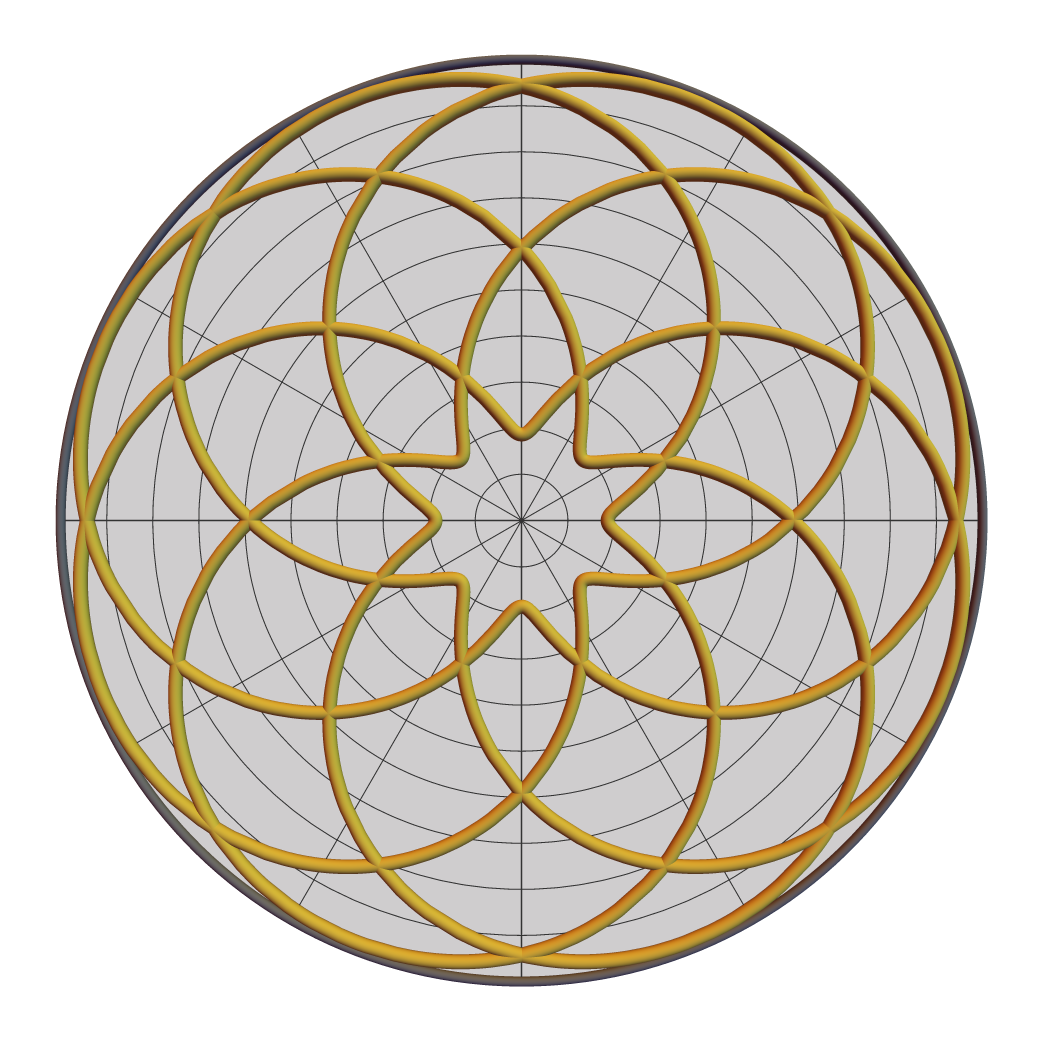}
	\includegraphics[height=5.55cm]{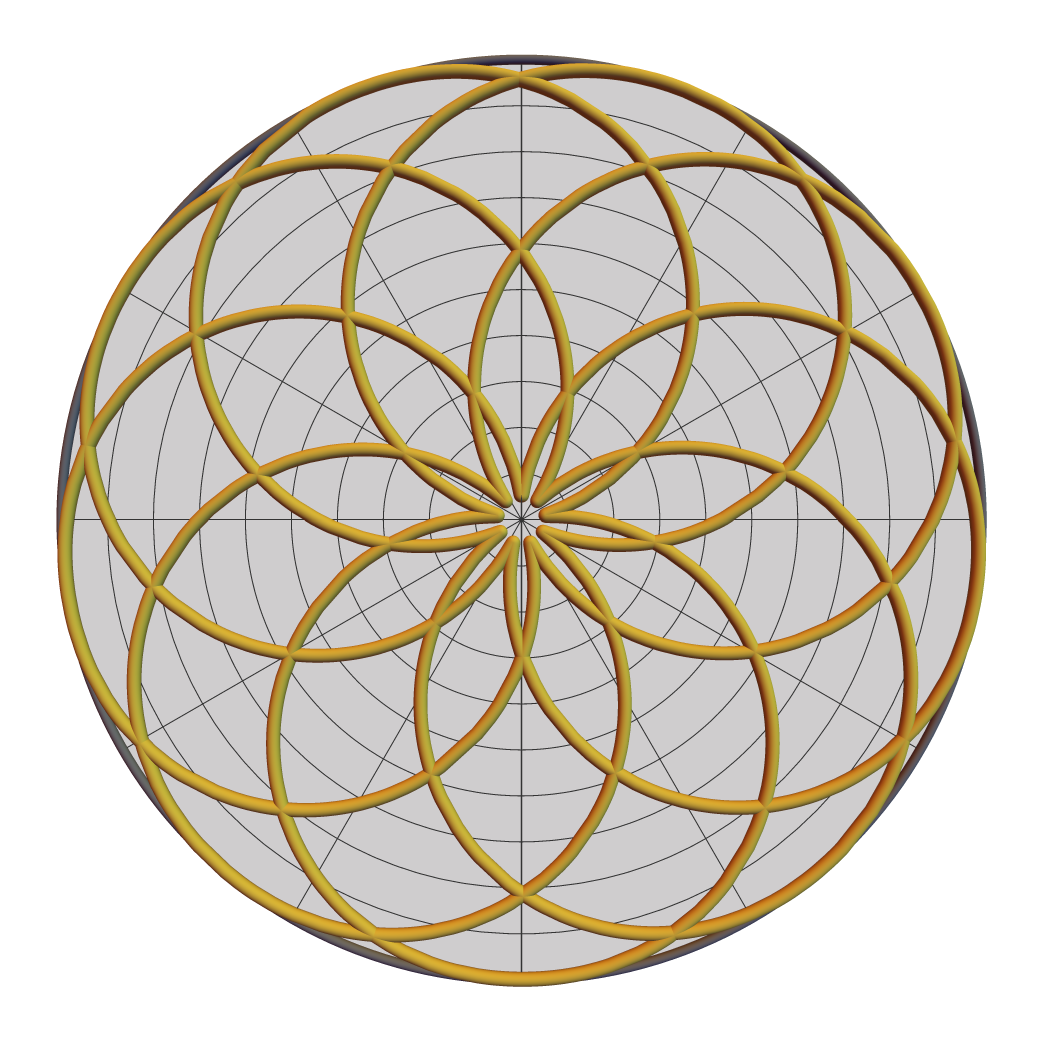}
	\includegraphics[height=5.55cm]{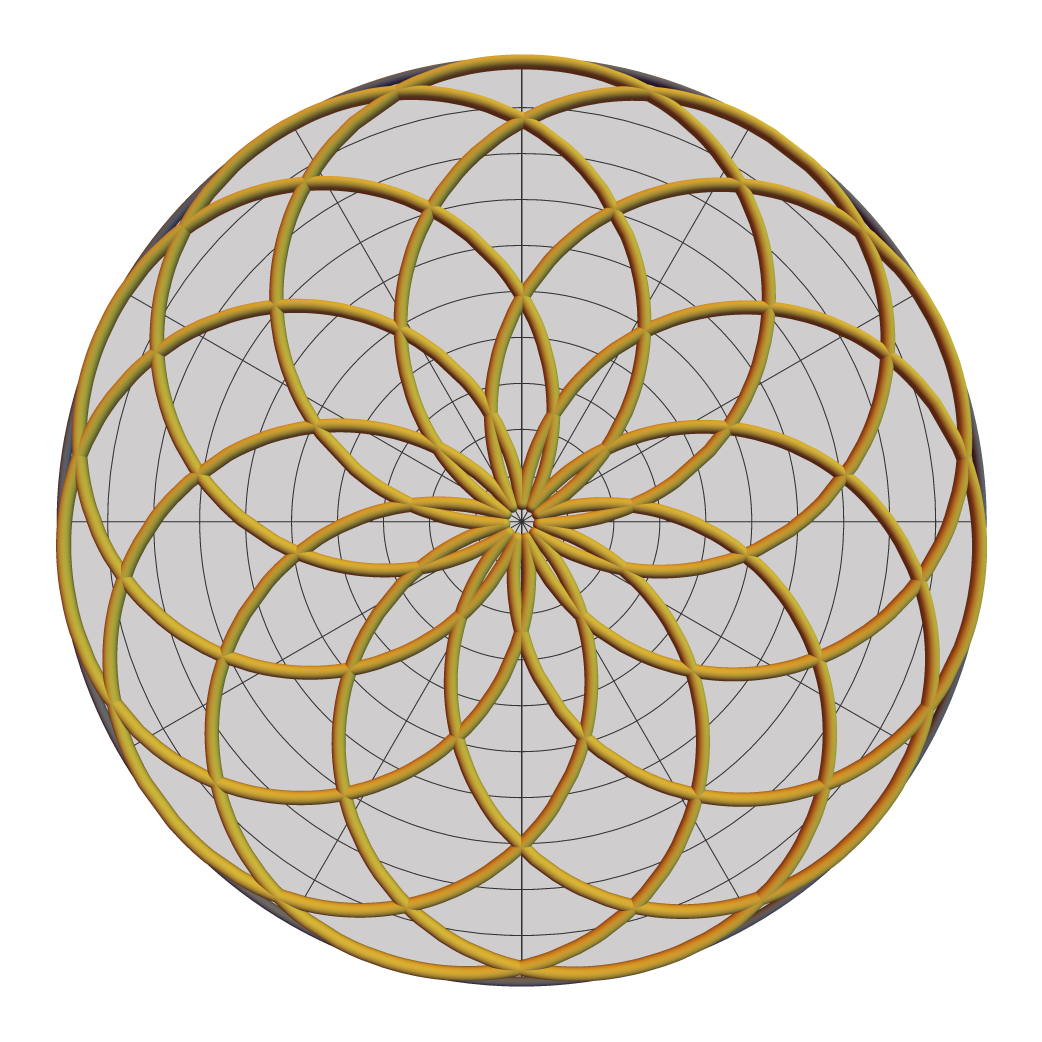}
	\caption{Closed (space-like) pseudo-hyperbolic $p$-elastic curves for $p=-1$ in $\mathbb{H}_1^2$ corresponding to the values $q=5/8$, $q=5/9$ and $q=6/11$, respectively. They are represented in the once-punctured unit disk $\mathring{\mathbb{D}}$. (Observe that the transformation that takes $\mathbb{H}_1^2\lvert_-$ to $\mathring{\mathbb{D}}$ involves an inversion.)}
	\label{F4}
\end{figure}

\begin{rem}\label{remH1} For the de Sitter $2$-space $\mathbb{H}_1^2$ (since $\epsilon=1$, then $p<0$ must hold) it is numerically evident that $\Lambda_p\to \pi$ when $a\to 0^-$, as is the case of the hyperbolic plane. Unfortunately, proving this assertion is more challenging since $\epsilon_2=(-1)^\epsilon=-1$ and, hence, $a+\epsilon_2p^2 u^{l-t}$ has no positive roots (c.f., proof of Theorem \ref{closed}). It is possible that a different choice of contour integrals would lead to the desired result. Nonetheless, the numerical evidence regarding the limit of $\Lambda_p$ when $a\to 0^-$ suggests that for every $p<0$ and pair of relatively prime natural numbers $(n,m)$ satisfying $m<2n<\sqrt{2}\,m$ there exists a closed $p$-elastic curve in $\mathbb{H}_1^2$ with nonconstant curvature (see Figures \ref{F2} and \ref{F4}).
\end{rem}

In Figures \ref{F1} and \ref{F2} we showed the closed $p$-elastic curves with nonconstant curvature $\gamma_{2,3}$, $\gamma_{3,5}$ and $\gamma_{4,7}$ in $\mathbb{H}_0^2$ for $p=3/2$ and in $\mathbb{H}_1^2$ for $p=-1$. These are the simplest choices for the parameters $n$ and $m$ satisfying $m<2n<\sqrt{2}\,m$. In Figure \ref{F3} we represent the closed $p$-elastic curves $\gamma_{5,8}$, $\gamma_{5,9}$ and $\gamma_{6,11}$ in $\mathbb{H}_0^2$ for $p=3/2$ (the corresponding curves in $\mathbb{H}_1^2$ for $p=-1$ are shown in Figure \ref{F4}).

\begin{figure}[h!]
	\centering
	\includegraphics[height=5.5cm]{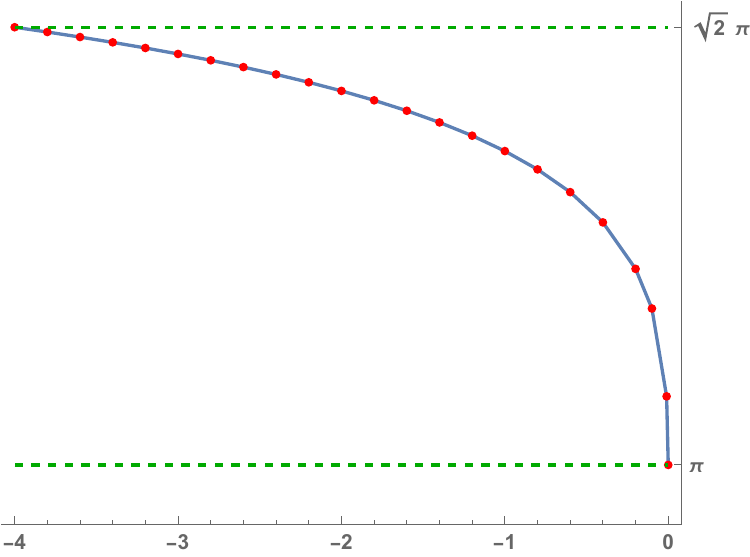}\quad\quad\quad\quad
	\includegraphics[height=5.5cm]{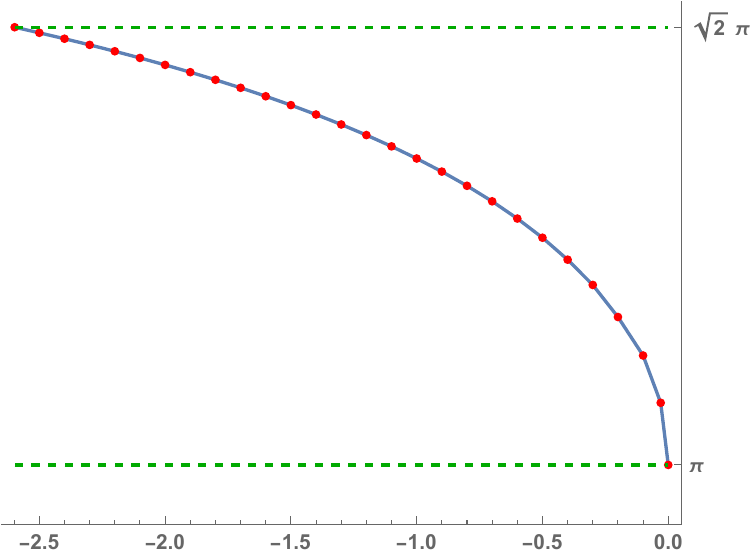}
	\caption{The numerical representation of the graph of the function $\Lambda_p:(a_*,0)\longrightarrow\mathbb{R}$ for $p=-1$ (Left) and $p=3/2$ (Right). For the case $p=-1$, $a_*=-4$; while for $p=3/2$, $a_*=-3\sqrt{3}/2$. We observe that in both cases the function $\Lambda_p$ decreases from $\sqrt{2}\,\pi$ (the limit when $a\to a_*$) to $\pi$ (the limit when $a\to 0$).}
	\label{FNumerical}
\end{figure}

\begin{rem}
	According to the numerical experiments the function $\Lambda_p:(a_*,0)\longrightarrow\mathbb{R}$ is monotonic (see Figure \ref{FNumerical} for the graph of $\Lambda_p$ in two particular cases). The validity of this assertion would imply that $\Lambda_p\left((a_*,0)\right)=(\pi,\sqrt{2}\pi)$ and for each $2\pi q\in(\pi,\sqrt{2}\pi)$, there exists a unique $a_q\in(a_*,0)$ such that $\Lambda_p(a_q)=2\pi q$. This would mean that for every $q\in\mathbb{Q}$ satisfying $1<2q<\sqrt{2}$, and under the restrictions on the parameter $p\in\mathbb{R}$ of Lemma \ref{lemma}, there exists a unique closed $p$-elastic curve in $\mathbb{H}_\epsilon^2$ with nonconstant curvature, up to isometries. In addition, none of them would be embedded.
\end{rem}

To prove the monotonicity of $\Lambda_p$ will require differentiating \eqref{Lambda} (see the expression \eqref{terrible} as well) with respect to $a$. The function $\Lambda_p$ represents, in general, a combination of extended hyperelliptic integrals where the coefficients depend on the parameter $a$ and, hence, differentiating $\Lambda_p$ directly seems prohibitively difficult. However, for the cases $p=2$ and $p=3$, the function $\Lambda_p$ can be written as a linear combination of standard complete elliptic integrals and so the derivative can be carried out employing the well known formulas for the derivatives of these elliptic integrals (see \cite{ABG0,LS} for these results). 

It is possible to find other values of $p\in(-\infty,0)\cup(1,\infty)$ for which the function $\Lambda_p$, possibly after a suitable change of variable, is an elliptic integral as well. Thus, the computation of its derivative would be feasible too. Although more complicated than the previous known cases ($p=2$ and $p=3$), one such a value is $p=3/2$. In the following result we show the monotonicity of $\Lambda_{3/2}$ by rewriting it in terms of standard complete elliptic integrals of the first and third kind. We first recall that the complete elliptic integrals of first $K$, second $E$, and third $\Pi$ kind are defined, respectively, by
\begin{eqnarray*}\label{KPi}
K(\zeta)&:=&\int_0^1\frac{du}{\sqrt{(1-u^2)(1-\zeta^2u^2)}}\,,\\
E(\zeta)&:=&\int_0^1\frac{\sqrt{1-\zeta^2u^2}}{\sqrt{1-u^2}}\,du\,,\\
\Pi(\chi,\zeta)&:=&\int_0^1 \frac{du}{(1-\chi u^2)\sqrt{(1-u^2)(1-\zeta^2u^2)}}\,.
\end{eqnarray*}
For later use, we also recall here that when $\zeta^2<\chi<1$ holds, the standard elliptic integral of the third kind $\Pi(\chi,\zeta)$ can be described involving the \emph{Heuman's Lambda} $\widehat{\Lambda}$ as
\begin{equation}\label{Heuman}
	\Pi(\chi,\zeta)=\frac{\pi}{2}\sqrt{\frac{\chi}{(1-\chi)(\chi-\zeta^2)}\,}\,\,\widehat{\Lambda}\left(\arcsin\sqrt{\frac{\chi-\zeta^2}{\chi(1-\zeta^2)}},\zeta\right).
\end{equation}
This expression follows after combining formulas $110.08$ and $150.01$ of \cite{BF}.

We are now in a position to prove the monotonicity of $\Lambda_{3/2}$.

\begin{theorem}\label{5.3} Let $p=3/2$ and so $\epsilon_2=1$. Then, the function $\Lambda_{3/2}:(a_*=-3\sqrt{3}/2,0)\longrightarrow\mathbb{R}$ defined in \eqref{Lambda} is given by
	\begin{equation}\label{simple}
		\Lambda_{3/2}(a)=\frac{2\sqrt{\alpha\beta(\alpha+\beta)}}{3\sqrt{2\alpha+\beta\,}}\,K(\zeta)+\pi\,\widehat{\Lambda}\left(\arcsin\sqrt{\frac{\chi-\zeta^2}{\chi(1-\zeta^2)}},\zeta\right).
	\end{equation}
where $\alpha>\beta>0$ are the only two positive roots of the polynomial $\widetilde{Q}_{3/2,a}(\kappa)=(-\kappa^3+9\kappa+4a)/4$ and
$$\zeta:=\sqrt{\frac{\alpha-\beta}{2\alpha+\beta}}\,,\quad\quad\quad \chi:=\frac{9(\alpha-\beta)}{9\alpha-\alpha\beta(\alpha+\beta)}\,.$$
Moreover, the function $\Lambda_{3/2}$ decreases monotonically from $\sqrt{2}\,\pi$ to $\pi$.
\end{theorem}
\textit{Proof.} Let $p=3/2$ (hence, $\epsilon_2=1$ necessarily) and consider the definition of the function $\Lambda_{p}$ given in \eqref{Lambda}. After substituting $p=3/2$ and using the first integral \eqref{fi} to make a change of variable (see also \eqref{terrible}), we obtain
\begin{eqnarray*}
\Lambda_{3/2}(a)&=&6\sqrt{-a}\int_\beta^\alpha \frac{\kappa\, d\kappa}{(4a+9\kappa)\sqrt{-\kappa^3+9\kappa+4a}}=\frac{2}{3}\sqrt{-a}\int_\beta^\alpha \frac{(4a+9\kappa)-4a}{(4a+9\kappa)\sqrt{-\kappa^3+9\kappa+4a}}\,d\kappa
\\&=&\frac{2}{3}\sqrt{-a}\left(\int_\beta^\alpha \frac{d\kappa}{\sqrt{-\kappa^3+9\kappa+4a}}-\frac{4a}{9}\int_\beta^\alpha \frac{d\kappa}{(\frac{4a}{9}+\kappa)\sqrt{-\kappa^3+9\kappa+4a}}\right).
\end{eqnarray*}
From formula 3.131-6 of \cite{GR} (see also formula 235.00 of \cite{BF}) the first integral in the second line above is a multiple of a suitable complete elliptic integral of the first kind, while from formula 3.137-6 of \cite{GR} (see also formula 236.02 of \cite{BF}) the second integral is a multiple of a suitable complete elliptic integral of the third kind. Combining both things, we deduce
\begin{equation}\label{aux}
	\Lambda_{3/2}(a)=\frac{2}{3}\sqrt{-a}\left(\frac{2}{\sqrt{\alpha-\delta}}\,K(\zeta)-\frac{8a}{(9\alpha+4a)\sqrt{\alpha-\delta}}\,\Pi(\chi,\zeta)\right),
\end{equation}
where $\zeta:=\sqrt{\alpha-\beta}/\sqrt{\alpha-\delta}$ and $\chi:=9(\alpha-\beta)/(9\alpha+4a)$. In both cases, $\delta<0<\beta<\alpha$ are the roots of the polynomial $-\kappa^3+9\kappa+4a$. Then, we have the Cardano-Vieta relations
$$\alpha+\beta+\delta=0\,,\quad\quad\quad\alpha\beta+\alpha\delta+\beta\delta=-9\,,\quad\quad\quad\alpha\beta\delta=4a\,.$$
Solving $\delta=-\alpha-\beta$ from the first relation and substituting this together with the last relation in \eqref{aux} we obtain the expression 
$$\Lambda_{3/2}(a)=\frac{2\sqrt{\alpha\beta(\alpha+\beta)}}{3\sqrt{2\alpha+\beta}}\left(K(\zeta)+\frac{\beta(\alpha+\beta)}{9-\beta(\alpha+\beta)}\,\Pi(\chi,\zeta)\right).$$

Moreover, we notice that above Cardano-Vieta relations can be solved in terms of just the largest root $\alpha$ as
$$a=\frac{1}{4}\alpha(\alpha^2-9)\,,\quad\quad\quad\beta=\frac{1}{2}\left(\sqrt{36-3\alpha^2}-\alpha\right),\quad\quad\quad\delta=-\frac{1}{2}\left(\sqrt{36-3\alpha^2}+\alpha\right).$$
Using this and \eqref{Heuman}, we can simplify the expression of $\Lambda_{3/2}$, obtaining \eqref{simple}.

We now observe that as $a\to a_*=-3\sqrt{3}/2$, then both $\alpha$ and $\beta$ tend to the curvature of the $p$-elastic circle given in \eqref{constantcurvature}. That is, for $p=3/2$, we have that $\alpha,\beta\to\sqrt{3}$ as $a\to -3\sqrt{3}/2$. Therefore, it follows that $\zeta\to 0$ and
$$\sqrt{\frac{\chi-\zeta^2}{\chi(1-\zeta^2)}}\longrightarrow \frac{2\sqrt{2}}{3}\,.$$
We then compute from \eqref{simple},
$$\lim_{a\to -3\sqrt{3}/2} \Lambda_{3/2}(a)=\frac{2\sqrt{6\sqrt{3}\,}}{3\sqrt{3\sqrt{3}\,}}\,K(0)+\pi\,\widetilde{\Lambda}\left(\arcsin \frac{2\sqrt{2}}{3},0\right)=\frac{\sqrt{2}}{3}\pi+\frac{2\sqrt{2}}{3}\pi=\sqrt{2}\,\pi\,,$$
where we have used the well known values $K(0)=\pi/2$ and $\widetilde{\Lambda}(\phi,0)=\sin\phi$ (see formulas 111.02 and 151.01 of \cite{BF}, respectively).

On the other hand, when $a\to 0$, then we deduce from the Cardano-Vieta relations that $\alpha\to 3$ while $\beta\to 0$. Therefore, $\zeta\to \sqrt{2}/2$ and
$$\sqrt{\frac{\chi-\zeta^2}{\chi(1-\zeta^2)}}\longrightarrow 1\,.$$
Observe that $K(\sqrt{2}/2)<\infty$ and that the coefficient of $K(\zeta)$ in \eqref{simple} tends to zero. Hence,
$$\lim_{a\to 0}\Lambda_{3/2}(a)=\pi\,\widetilde{\Lambda}\left(\frac{\pi}{2},\frac{\sqrt{2}}{2}\right)=\pi\,,$$
where we have used that $\widetilde{\Lambda}(\pi/2,\zeta)=1$ (see formula 151.01 of \cite{BF}).

It remains to prove the monotonicity of the function $\Lambda_{3/2}(a)$. To this end, we will compute the derivative $\Lambda_{3/2}'(a)$ of \eqref{simple} with respect to $a$ and check its sign. We first observe that the largest root $\alpha$ of the polynomial $-\kappa^3+9\kappa+4a$ increases monotonically from $\sqrt{3}$ to $3$ as $a$ varies in $(-3\sqrt{3}/2,0)$. Indeed, implicitly differentiating $a=\alpha(\alpha^2-9)/4$, we get
$$\alpha'(a)=\frac{4}{3(\alpha^2-3)}>0\,.$$
Therefore, the sign of $\Lambda_{3/2}'(a)$ is the same as that of the derivative of $\Lambda_{3/2}$, \eqref{simple}, with respect to $\alpha$, which is simpler to compute because all the functions involved can be described in terms of $\alpha$ in a simpler way (see the above Cardano-Vieta relations).

Employing the formulas 710.11 and 730.04 of \cite{BF} to compute the derivative of the Heuman's Lambda $\widehat{\Lambda}$ and formula 710.00 of \cite{BF} to compute the derivative of $K(\zeta)$, it is possible to obtain the derivative of $\Lambda_{3/2}$ with respect to $\alpha$ just as a linear combination of the complete elliptic integrals of first and second kind, $K(\zeta)$ and $E(\zeta)$, respectively. In addition, long simplifications of this expression by means of the Cardano-Vieta relations allow us to check that the derivative of \eqref{simple} with respect to $\alpha$ is everywhere negative when $\alpha\in(\sqrt{3},3)$. (Even after these simplifications and although this derivative only involves two terms, namely, a multiple of $K(\zeta)$ and a multiple of $E(\zeta)$, the expression still very long and will necessitate of several lines to describe it. Hence, we have opted here to omit it.) The fact that the derivative of \eqref{simple} with respect to $\alpha$ is negative means that the sign of $\Lambda_{3/2}'(a)$ is also negative and, hence, $\Lambda_{3/2}(a)$ decreases monotonically from $\sqrt{2}\,\pi$ to $\pi$, which concludes the proof. \hfill$\square$

\subsection{Energy}

We next study the energy of closed $p$-elastic curves in $\mathbb{H}_\epsilon^2$ and show, in two particular cases, how this value can be computed employing standard elliptic integrals. 

Let $\gamma_a\equiv\gamma_{n,m}$, $0>a>a_*$, be a closed $p$-elastic curve closing up in $m$ periods of the curvature. Then, using \eqref{fi} to make a change of variable in the definition of the $p$-elastic functional $\mathbf{\Theta}_p$, we have
\begin{equation}\label{energyhyperbolic}
	\mathbf{\Theta}_p(\gamma_a)=2\,m\,p\left(p-1\right)\int_\beta^\alpha\frac{\kappa^{2(p-1)}}{\sqrt{a-(1-p)^2\kappa^{2p}+p^2\kappa^{2(p-1)}}}\,d\kappa\,,
\end{equation}
if $\gamma_a$ is a hyperbolic curve (i.e., $\epsilon=0$) in which case $p>1$, and
\begin{equation}\label{energypseudohyperbolic}
	\mathbf{\Theta}_p(\gamma_a)=2\,m\,p\left(p-1\right)\int_\beta^\alpha\frac{d\kappa}{\kappa^{1-p}\sqrt{a\kappa^{2(1-p)}+(p-1)^2\kappa^2-p^2\,}}\,,
\end{equation}
if $\gamma_a$ is a pseudo-hyperbolic curve (i.e., $\epsilon=1$) for which $p<0$ holds.

In both cases, as already done in the proof of Lemma \ref{lemma} employing Corollary 4.2 of \cite{Pe}, we can compute the limit when $a\to a_*^+$, obtaining after straightforward simplifications
\begin{equation}\label{limit}
	\lim_{a\to a_*^+}\mathbf{\Theta}_p(\gamma_a)=\sqrt{2}\,m\left((-1)^\epsilon p\right)^{p/2}\left((-1)^\epsilon \left(p-1\right)\right)^{(1-p)/2}\pi=\sqrt{-2a_*}\,m\,\pi\,.
\end{equation}

Observe that for suitable choices of the energy parameter $p\in(-\infty,0)\cup(1,\infty)$ the integrals arising in \eqref{energyhyperbolic} and \eqref{energypseudohyperbolic}, respectively, can be described in terms of standard complete elliptic integrals. We illustrate this fact for two particular cases, namely $p=-1$ and $p=3/2$, in the following result. 

\begin{proposition} Let $\gamma_a\equiv\gamma_{m,n}$ be a closed $p$-elastic curve in $\mathbb{H}^2_\epsilon$ closing up in $m$ periods of its curvature $\kappa$ and denote by $\alpha>\beta>0$ the maximum and minimum values of $\kappa$. Then, the value of the $p$-elastic functional $\mathbf{\Theta}_p$ is given by:
	\begin{enumerate}[(i)]
		\item If $p=-1$,
		\begin{equation}\label{Theta-1}
			\mathbf{\Theta}_{-1}(\gamma_a)=4m\,\frac{\beta}{\alpha^2}\,\Pi(\zeta^2,\zeta)\,,
		\end{equation}
		where $\zeta=\sqrt{\alpha^2-\beta^2}/\alpha$.
		\item If $p=3/2$,
		\begin{equation}\label{Theta3/2}
			\mathbf{\Theta}_{3/2}(\gamma_a)=6m\left(\sqrt{2\alpha+\beta} \,E(\zeta)-\frac{\alpha+\beta}{\sqrt{2\alpha+\beta}} \,K(\zeta)\right),
		\end{equation}
		where $\zeta=\sqrt{\alpha-\beta}/\sqrt{2\alpha+\beta}$.
	\end{enumerate}
\end{proposition}
\textit{Proof.} Consider first the case $p=-1$. Substituting this on the expression \eqref{energypseudohyperbolic} for the general $p$-elastic functional $\mathbf{\Theta}_p$, $p<0$, we have
$$\mathbf{\Theta}_{-1}(\gamma_a)=4m\int_\beta^\alpha \frac{d\kappa}{\kappa^2\sqrt{a\kappa^4+4\kappa^2-1}}=\frac{2m}{\sqrt{-a}}\int_{\beta^2}^{\alpha^2}\frac{du}{u\sqrt{(\alpha^2-u)(u-\beta^2)u}}\,,$$
where in the last equality we have applied the change of variable $u=\kappa^2$. Observe that since $\alpha>\beta>0$ are the only two positive roots of the polynomial $a\kappa^4+4\kappa^2-1$, then $\alpha^2>\beta^2>0$ are the only two roots of $au^2+4u-1$. It then follows from formula 3.137-6 of \cite{GR} (see also formula 236.02 of \cite{BF}), that
$$\mathbf{\Theta}_{-1}(\gamma_a)=\frac{4m}{\alpha^3\sqrt{-a}}\,\Pi(\zeta^2,\zeta)\,,$$
where $\zeta=\sqrt{\alpha^2-\beta^2}/\alpha$. Finally, from the Cardano-Vieta relation $\alpha^2\beta^2=-1/a$ we deduce that
$$\frac{4m}{\alpha^3\sqrt{-a}}=4m\frac{\beta}{\alpha^2}\,,$$
and prove the statement.

Assume now that $p=3/2$ holds. From the expression of the general $p$-elastic functional $\mathbf{\Theta}_p$, $p>1$, given in \eqref{energyhyperbolic} and formula 3.132-5 of \cite{GR} (see also formula 236.16 of \cite{BF}), we have that
\begin{eqnarray*}
	\mathbf{\Theta}_{3/2}(\gamma_a)&=&3m\int_\beta^\alpha \frac{\kappa}{\sqrt{-\kappa^3+9\kappa+4a}}\,d\kappa=3m\int_\beta^\alpha \frac{\kappa}{\sqrt{(\kappa-\delta)(\kappa-\beta)(\alpha-\kappa)}}\,d\kappa\\
	&=&6m\left(\frac{\delta}{\sqrt{\alpha-\delta}}\,K(\zeta)+\sqrt{\alpha-\delta}\,E(\zeta)\right),
\end{eqnarray*}
where $\zeta=\sqrt{\alpha-\beta}/\sqrt{\alpha-\delta}$ and $\alpha>\beta>0>\delta$ are the roots of the polynomial $-\kappa^3+9\kappa+4a$. We next use the Cardano-Vieta relations (see the proof of Theorem \ref{5.3}) to deduce that $\delta=-\alpha-\beta$ and, hence, conclude the proof. \hfill$\square$

\begin{rem} From the analysis of the proof of Theorem \ref{periodic}, we deduce that when $a\to a_*^+$ both $\alpha$ and $\beta$ tend to the curvature of a non-geodesic $p$-elastic circle given in \eqref{constantcurvature}. Therefore, from the expressions \eqref{Theta-1} and \eqref{Theta3/2}, and the well known values $K(0)=E(0)=\Pi(0,0)=\pi/2$ (see formulas 111.02 and 111.06 of \cite{BF}) we conclude with
	\begin{eqnarray*}
		\lim_{a\to a_*^+} \mathbf{\Theta}_{-1}(\gamma_a)&=&2\sqrt{2\,}\,m\,\pi\,,\\
		\lim_{a\to a_*^+}\mathbf{\Theta}_{3/2}(\gamma_a)&=&\sqrt{3\sqrt{3}\,}\,m\,\pi\,,
	\end{eqnarray*}
which coincides with the general limit computed in \eqref{limit}.
\end{rem}

\subsection{Kinematics}

In this part we describe an interesting behavior of the closed $p$-elastic curves in $\mathbb{H}_\epsilon^2$ analyzed in this paper as the energy parameter $p$ varies in its respective domain ($p>1$ if $\epsilon=0$ and $p<0$ if $\epsilon=1$). According to the numerical experiments about the monotonicity of $\Lambda_p$, the closed $p$-elastic curves of this paper are all.

These closed $p$-elastic curves in $\mathbb{H}_\epsilon^2$ are either circles (Proposition \ref{constant}) or they arise from a pair of relatively prime natural numbers $(n,m)$ satisfying $m<2n<\sqrt{2}\,m$ (Corollary \ref{cor} and Remark \ref{remH1}). Observe that the latter condition is independent of the parameter $p\in\mathbb{R}$. Hence, one can fix the pair $(n,m)$ and consider the evolution of the associated family of closed $p$-elastic curves as $p$ varies.

To clarify this evolution, we will consider first the case of $p$-elastic circles. From Proposition \ref{constant}, we obtain that the radii of $p$-elastic circles in $\mathbb{H}_\epsilon^2$, viewed as curves in $\mathbb{L}^3$, is $r\equiv r(p)=\sqrt{(-1)^\epsilon(p-1)}$. The function $r(p)$ is monotonic in $p$ for each possible choice of ambient space $\mathbb{H}_\epsilon^2$, $\epsilon=0,1$. Indeed:
\begin{enumerate}[(i)]
	\item If $\epsilon=0$, the radii of the $p$-elastic circles in the hyperbolic plane $\mathbb{H}^2\subset\mathbb{L}^3$ increases monotonically from $0$ to $\infty$, as $p$ varies in $(1,\infty)$. This means that $p$-elastic circles are, up to isometries, parallel circles in $\mathbb{H}^2$ whose height increases and they vary from the pole $(0,0,1)\in\mathbb{R}^3$ to the boundary at infinity. In the Poincar\'e disk model, $p$-elastic circles are circles centered at the origin which go expanding from $(0,0)$ to the unit circle, i.e., the ideal boundary of $\mathbb{D}$.
	\item If $\epsilon=1$, the radii of the $p$-elastic circles in the de Sitter $2$-space $\mathbb{H}_1^2\subset\mathbb{L}^3$ decreases monotonically from $\infty$ to $1$, as $p$ varies in $(-\infty,0)$. This means that $p$-elastic circles are, up to isometries, parallel circles in $\mathbb{H}_1^2$ whose height increases and they vary from the boundary at infinity of the bottom half of $\mathbb{H}_1^2$ to the equator $\mathbb{H}_1^2\cap\{z=0\}$. In the model $\mathring{\mathbb{D}}$ for the bottom half of $\mathbb{H}_1^2$, $p$-elastic circles are circles centered at the origin which go expanding from the limiting point $(0,0)$ to the unit circle, i.e., the boundary of $\mathring{\mathbb{D}}$. (Observe that the transformation that takes $\mathbb{H}_1^2\lvert_-$ to $\mathring{\mathbb{D}}$ involves an inversion.)
\end{enumerate}

The evolution of closed $p$-elastic curves with nonconstant curvature imitates that of $p$-elastic circles. For a fixed pair of relatively prime natural numbers $(n,m)$ satisfying $m<2n<\sqrt{2}\,m$, consider the family of closed $p$-elastic curves in $\mathbb{H}_\epsilon^2$ with nonconstant curvature $\{\gamma_{n,m}\}_p$, where $p>1$ if $\epsilon=0$, or $p<0$ if $\epsilon=1$. The choice $(n,m)$ fixes both the winding number and the symmetry of all the curves $\gamma_{n,m}$ in the family, independently of $p$. In other words, in the disk models, all the curves in the family $\{\gamma_{n,m}\}_p$ go around the center of the disk $n$ times and have $m$ lobes.

\begin{figure}[h!]
	\makebox[\textwidth][c]{\centering
		\includegraphics[height=4.5cm]{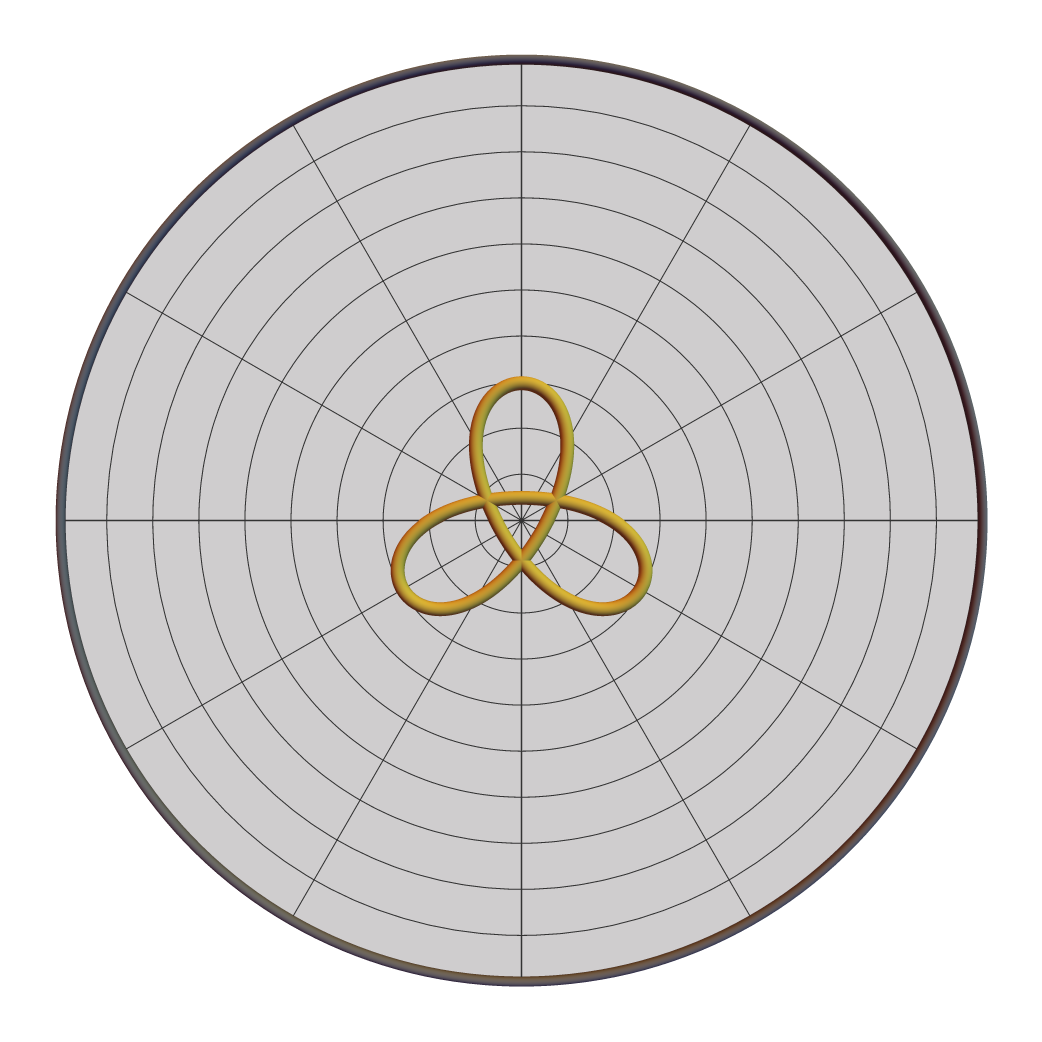}
		\includegraphics[height=4.5cm]{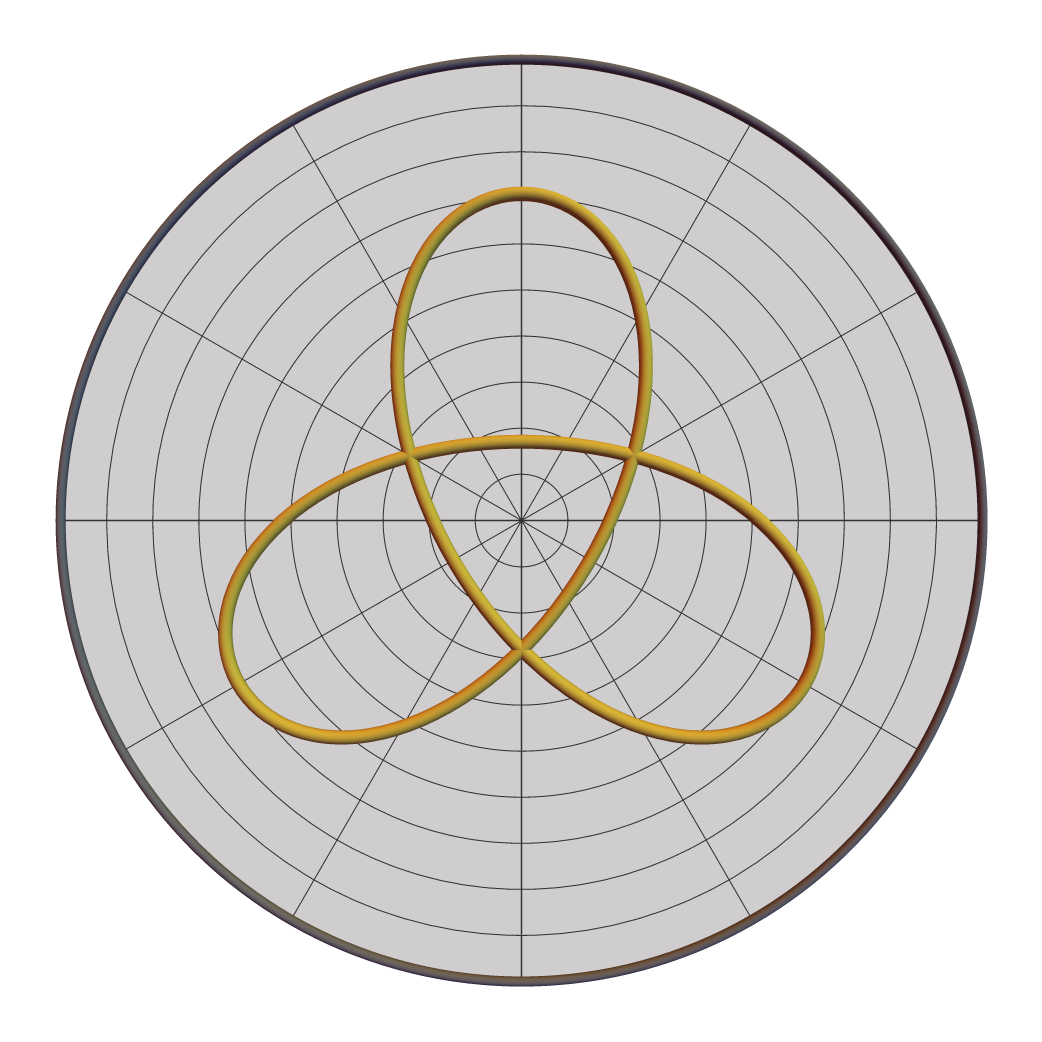}
		\includegraphics[height=4.5cm]{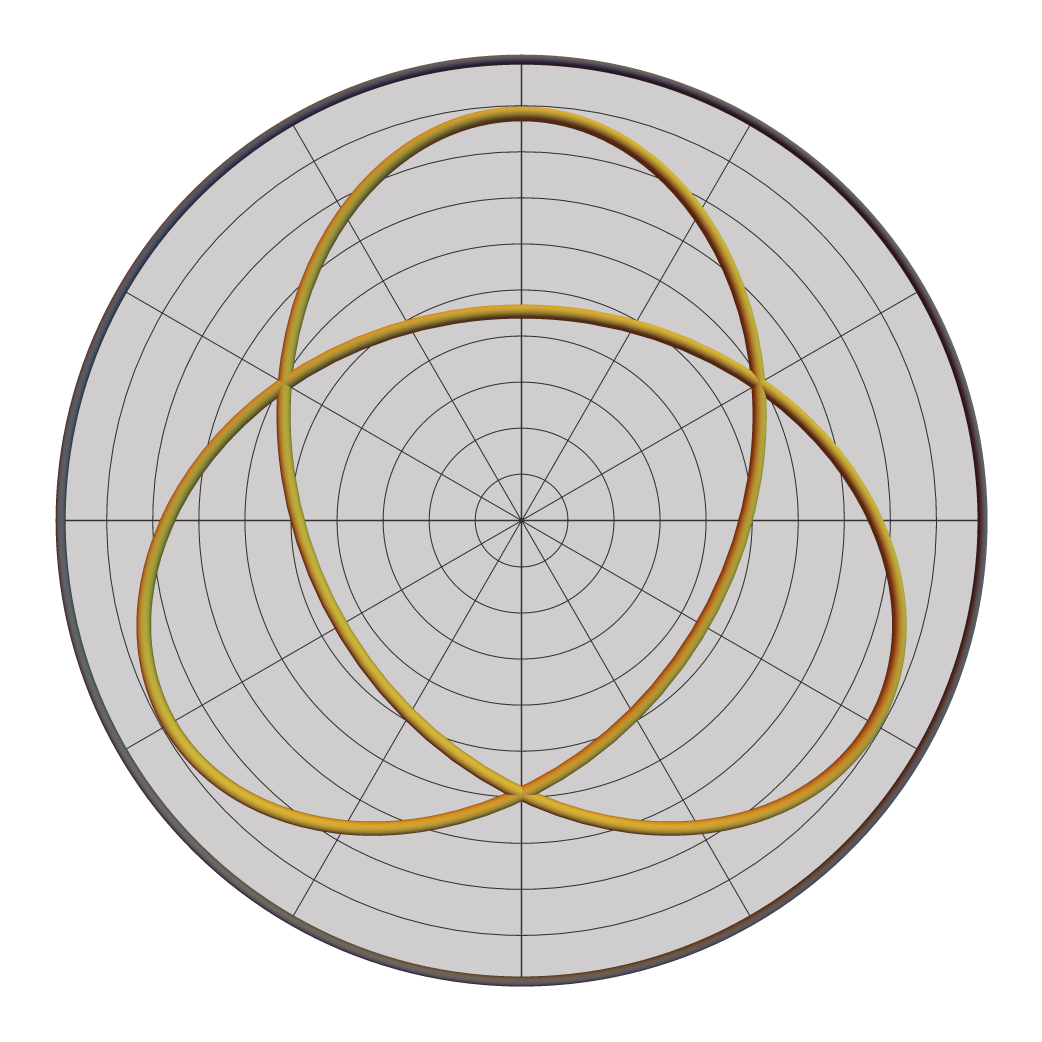}
		\includegraphics[height=4.5cm]{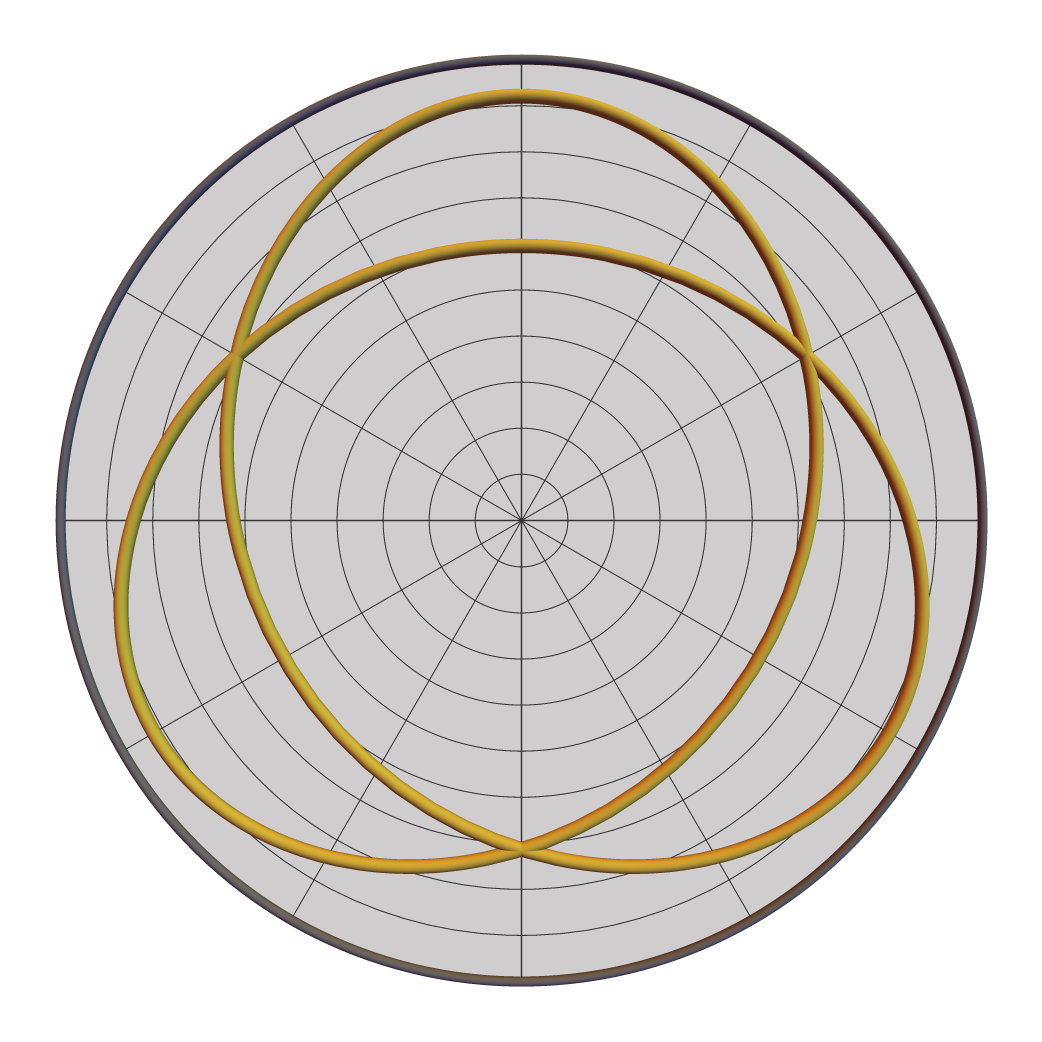}}
	\caption{Closed hyperbolic $p$-elastic curves in $\mathbb{H}_0^2$ of type $\gamma_{2,3}$ for different values of $p>1$. From left to right: $p=1.1$, $p=2$, $p=7$ and $p=15$.}
	\label{Evolution}
\end{figure}

Assume first that $\epsilon=0$ holds, then in the Poincar\'e disk model of the hyperbolic plane $\mathbb{H}^2$ the closed $p$-elastic curves $\gamma_{n,m}$ go expanding as $p$ varies in $(1,\infty)$, from the center $(0,0)$ (limiting case when $p\to 1^+$) to the unit circle (when $p\to\infty$). In Figure \ref{Evolution} we show four closed hyperbolic $p$-elastic curves for different values of $p$ and for the same type, namely, for $n=2$ and $m=3$. This figure illustrates the evolution described above.

\begin{figure}[h!]
	\makebox[\textwidth][c]{\centering
		\includegraphics[height=4.5cm]{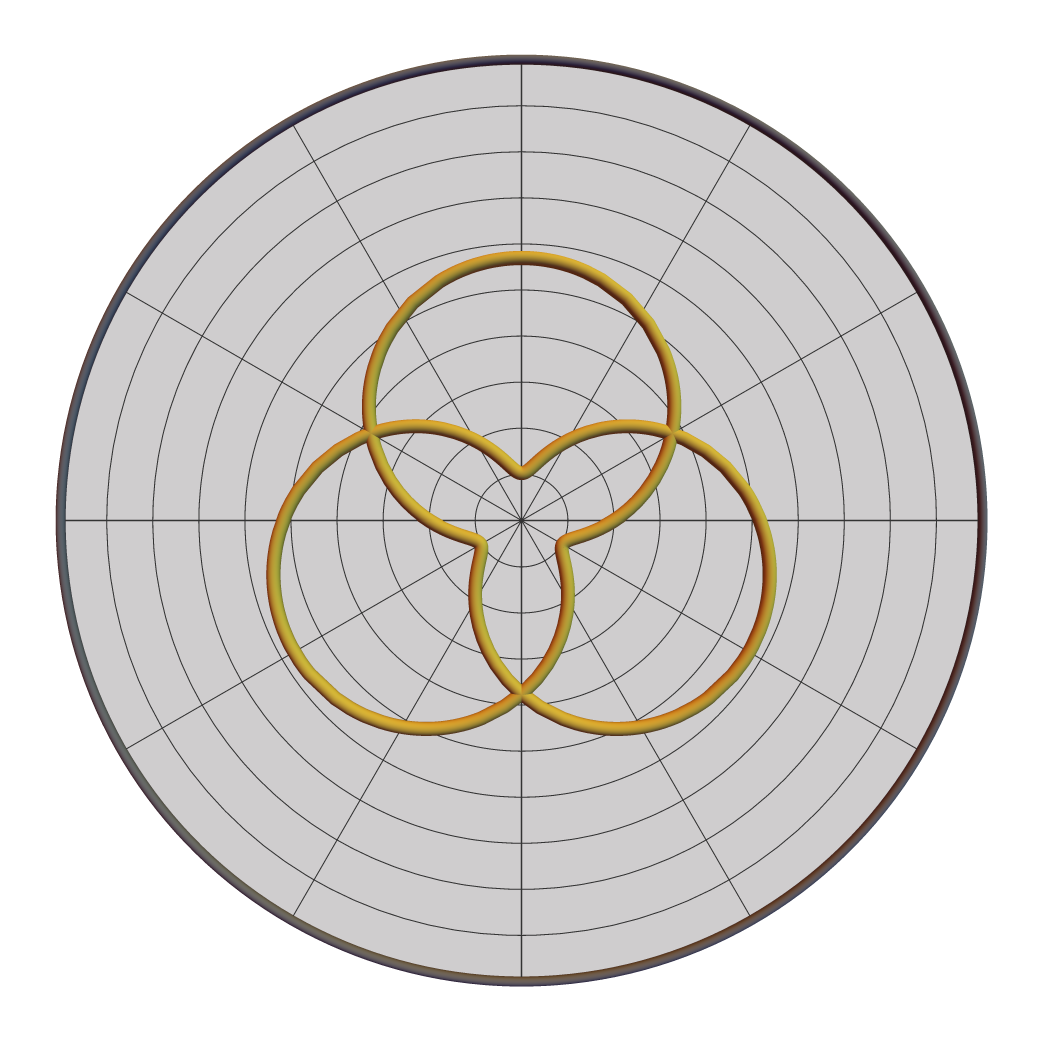}
		\includegraphics[height=4.5cm]{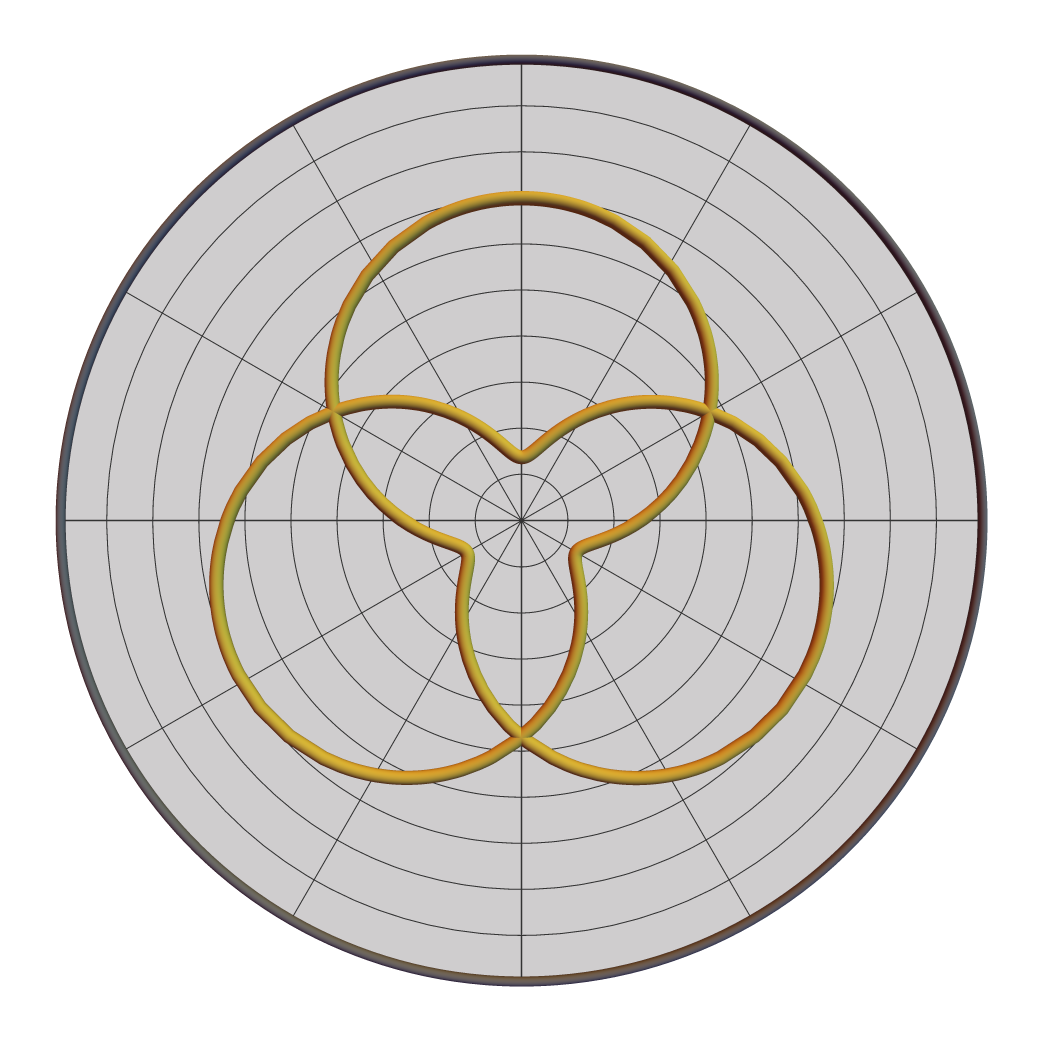}
		\includegraphics[height=4.5cm]{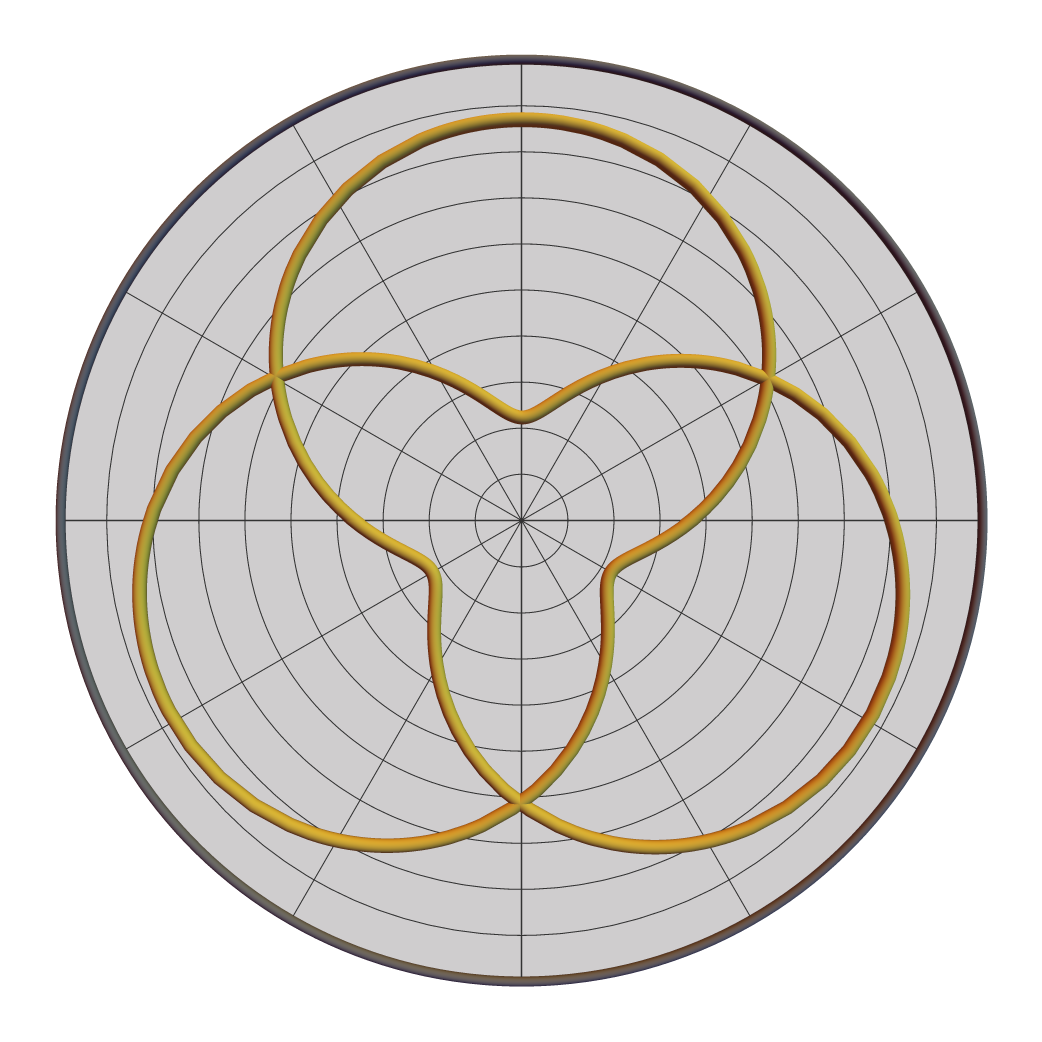}
		\includegraphics[height=4.5cm]{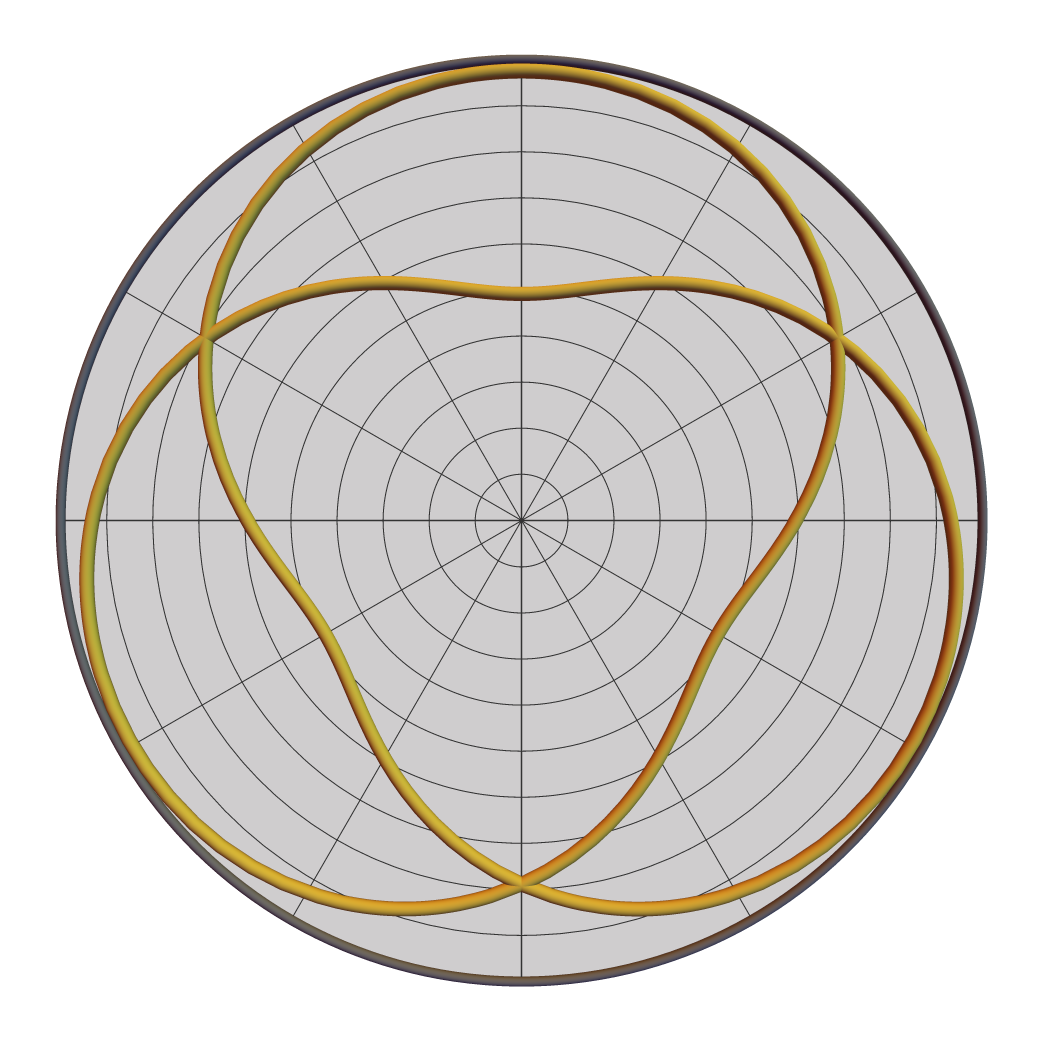}}
	\caption{Closed (space-like) pseudo-hyperbolic $p$-elastic curves in $\mathbb{H}_1^2$ of type $\gamma_{2,3}$ for different values of $p<0$. From left to right: $p=-9$, $p=-5$, $p=-2$ and $p=-1/2$. (Observe that the transformation that takes $\mathbb{H}_1^2\lvert_-$ to $\mathring{\mathbb{D}}$ involves an inversion.)}
	\label{EvolutionP}
\end{figure}

On the other hand, if $\epsilon=1$, the closed $p$-elastic curves $\gamma_{n,m}$ in the model disk $\mathring{\mathbb{D}}$ for the bottom half of the de Sitter $2$-space $\mathbb{H}_1^2$ go expanding as $p$ varies in $(-\infty,0)$, from the center $(0,0)$ (when $p\to-\infty$) to the unit circle (when $p\to 0^-$). Figure \ref{EvolutionP} illustrates four closed pseudo-hyperbolic $p$-elastic curves in $\mathbb{H}_1^2$ of type $\gamma_{2,3}$ for different values of $p\in(-\infty,0)$.

\begin{rem} It is interesting to understand above described evolution when $\mathbb{H}_\epsilon^2$ is a subset of $\mathbb{R}^3$, that is, in the quadric models. Moreover, it will also give insight to combine above observations with those of \cite{GPT} about closed spherical $p$-elastic curves for $p\in(0,1)$. Indeed, this will give the continuous deformation for every $p\in\mathbb{R}$ (of course, with the remarkable exception of $p=1$, which is a ``singularity'') described in the last part of the Introduction.
\end{rem}

\section*{Acknowledgments}

The authors would like to thank Professor Alexander Solynin for his valuable comments and suggestions regarding the application of Cauchy's Integral Theorem and the choice of the suitable contour integrals. H. T's research is partially supported by grants from the National Science Foundation [DMS-2104988] and the Vietnam Institute for Advanced Study in Mathematics.

\begin{flushleft}
	\'Alvaro P{\footnotesize \'AMPANO}\\
	Department of Mathematics and Statistics, Texas Tech University, Lubbock, TX, 79409, USA\\
	E-mail: alvaro.pampano@ttu.edu
\end{flushleft}

\begin{flushleft}
	Miraj S{\footnotesize AMARAKKODY}\\
	Department of Mathematics and Statistics, Texas Tech University, Lubbock, TX, 79409, USA\\
	E-mail: miraj.samarakkody@ttu.edu
\end{flushleft}

\begin{flushleft}
	Hung T{\footnotesize RAN}\\
	Department of Mathematics and Statistics, Texas Tech University, Lubbock, TX, 79409, USA\\
	E-mail: hung.tran@ttu.edu
\end{flushleft}

\end{document}